\documentclass[12pt]{article}
\usepackage[american]{babel}
\usepackage{amsmath}
\usepackage{latexsym}
\usepackage{amssymb}
\usepackage{theorem}
\usepackage[arrow,matrix,curve]{xy}
\usepackage{diagrams}
\theoremstyle{change}
\newtheorem{Thm}{Theorem}[section]
\newtheorem{Cor}[Thm]{Corollary}
\newtheorem{Prop}[Thm]{Proposition}
\newtheorem{Lem}[Thm]{Lemma}
{\theorembodyfont{\rmfamily}
\newtheorem{Num}[Thm]{}

\newtheorem{Def}[Thm]{Definition}}
\newcommand{\proof}{\par\medskip\rm\emph{Proof. }}
\newcommand{\qed}{\ \hglue 0pt plus 1filll $\Box$}

\newcommand{\too}{\longrightarrow}
\newcommand{\mapstoo}{\longmapsto}

\renewcommand{\SS}{\mathbb{S}}
\newcommand{\RR}{\mathbb{R}}
\newcommand{\ZZ}{\mathbb{Z}}
\newcommand{\FF}{\mathbb{F}}
\newcommand{\HH}{\mathbb{H}}
\newcommand{\OO}{\mathbb{O}}
\newcommand{\CC}{\mathbb{C}}
\newcommand{\id}{\mathrm{id}}
\newcommand{\SKIP}[1]{}
\newcommand{\SO}{\mathrm{SO}}
\newcommand{\SU}{\mathrm{SU}}
\newcommand{\SL}{\mathrm{SL}}
\newcommand{\U}{\mathrm{U}}
\newcommand{\GL}{\mathrm{GL}}
\newcommand{\Sp}{\mathrm{Sp}}
\newcommand{\Spin}{\mathrm{Spin}}
\newcommand{\fsl}{\mathfrak{sl}}
\newcommand{\fso}{\mathfrak{so}}
\newcommand{\fsu}{\mathfrak{su}}
\newcommand{\fu}{\mathfrak{u}}
\newcommand{\fsp}{\mathfrak{sp}}
\newcommand{\fp}{\mathfrak{p}}
\newcommand{\fh}{\mathfrak{h}}
\newcommand{\fg}{\mathfrak{g}}
\newcommand{\fk}{\mathfrak{k}}
\newcommand{\fspin}{\mathfrak{spin}}
\newcommand{\I}{\mathtt{I}}
\newcommand{\ti}{\mathtt{i}}
\newcommand{\eps}{\varepsilon}
\renewcommand{\emptyset}{\varnothing}
\newcommand{\cP}{\mathcal P}
\newcommand{\cL}{\mathcal L}

\begin{document}

\title{\bf Two-transitive Lie groups}
\author{Linus Kramer\thanks{%
Supported by a Heisenberg fellowship by the Deutsche Forschungsgemeinschaft}} 
\maketitle

\begin{abstract}
{\em Using a characterization of parabolics in reductive
Lie groups due to Furstenberg, elementary
properties of buildings, and some algebraic topology,
we give a new proof of Tits' classification
of 2-transitive Lie groups.}
\end{abstract}
Among many other results, Tits classified in \cite{Ti}
all 2-transitive Lie groups.
His proof is based on Dynkin's classification of maximal
complex subalgebras in complex simple Lie algebras; it
is long and depends on consideration of various cases.
Since the resulting list of groups is also long (at least in the
affine case), it is clear that there cannot be a very short proof
of the full classification. On the other hand, Lie theory has
developed since the time \cite{Ti} was written. In particular,
Tits himself changed the picture through his theory of
buildings (as he pointed out, his paper \cite{Ti} was one of the
motivations for him to invent buildings). The language, the methods,
and the terminology have changed since then, and it is natural
to look for a new (and shorter) proof of Tits' 
classification. 
Note also that the proof presented in \cite{Ti} IV F 1.2, p.~222,
does not cover certain real forms of exceptional
groups --- a footnote on p.~223 asserts that
Tits found a proof for these cases, too, after the manuscript
went into print; see also \emph{loc.cit.} p.~240.
The details were never published.

Almost at the same time as Tits, Borel \cite{Bo}
determined all 1-connected spaces $X$ which admit a 2- or
3-transitive Lie group action; however, Borel did not classify
the corresponding groups. His proof relies on spectral
sequences, Freudenthal's theory of ends, and on the results
of Borel and De~Siebenthal about homogeneous spaces of positive
Euler characteristic.

In this paper, we give a complete proof for Tits' classification.
The main ingredients are a characterization
of parabolics in Lie groups due to Furstenberg,
elementary properties of buildings,
some algebraic topology (certainly more elementary than
the machinery employed in Borel's work \cite{Bo}), and representation
theory of semisimple (compact) Lie groups.

As we remarked before, the only published proof for the classification
is \cite{Ti}. The classification in the affine
case is also stated (but not proved) in V\"olklein \cite{Volk},
and some remarks on the
strategy of Tits' original proof can be found in Salzmann \emph{et al.}
\cite{CPP} 96.15 and 96.16.

Related results for other classes of
groups are
Knop's classification \cite{Knop} of 2-transitive actions of
algebraic groups over algebraically closed fields in arbitrary
characteristic (which is achieved by quite different methods), and
the classification of all finite 2-transitive
groups; see Dixon-Mortimer \cite{DM} 7.7 for a
description of these groups. In the course of our classification,
we recover Knop's result for the special case of complex
algebraic groups.

The main results of the classification are as follows.

\medskip\noindent
\textbf{Theorem A}
{\em
Let $G$ be a locally compact, $\sigma$-compact topological
transformation group acting effectively and 2-transitively on
a space $X$ which is not totally disconnected. Then
$G$ is a Lie group and $X\cong G/G_x$ is a connected manifold.
The connected component $G^\circ$ is simple if and only if
$X$ is compact.}

\smallskip\noindent
Note that it is not enough to assume that the group $G$ is
locally compact; the full homeomorphism group of any topological
manifold, endowed with the discrete topology, satisfies all the
other conditions of the theorem.

\medskip\noindent
\textbf{Theorem B}
{\em
If $(G,X)$ is as in Theorem A, and if $X$ is compact, then
$X$ is either the point set
of a projective space,
or the set of all absolute points of a polarity
(of index 1) in a projective space.
In the first case, $G$ is (a finite extension of)
the little projective group, and in the second case,
$G$ is (a finite extension of) the centralizer of the polarity in the
little projective group of the projective space.}

\smallskip\noindent
The projective spaces in question real, complex, quaternionic,
or octonionic, and the
possibilities for the groups $G$ are explicitly determined,
see Theorem \ref{CompactClass}.

\medskip\noindent
\textbf{Theorem C}
{\em
If $(G,X)$ is as in Theorem A, and if $X$ is noncompact, then
$X\cong\RR^m$ is a real vector space, and $G$ is a
semidirect product $G=G_x\ltimes\RR^m$, where $G_x\leq\GL_m\RR$
is a linear group acting transitively on the nonzero vectors.}

\smallskip\noindent
We determine explicitly the connected linear groups which act
transitively on the nonzero vectors in Theorem \ref{TrsLinGrp}.

\medskip
The dichotomy that either
$X$ is compact and $G^\circ$ is simple, or that $X$ is noncompact
and $G$ of affine type is proved by ideas similar to Borel's, but
with a modest amount of algebraic topology. The
classification of transitive linear groups depends very much on
representation theory.

In case where $X$ is compact, the key ingredient is a characterization
of parabolics due to Furstenberg. This characterization was
used by Burns-Spatzier \cite{BS} in order to classify compact 
connected buildings with strongly transitive automorphism
groups. Using Furstenberg's result, the classification
is reduced to a problem about the $W$-valued distance in spherical
buildings. Here, we need some elementary properties of buildings,
(but the proof does not depend on the classification of 
spherical Moufang buildings).
We rely of course on the classification of real simple Lie groups
and their structure theory.


\medskip
\noindent
\textbf{Outline of the classification.}

The first section collects some basic material on 2-transitive
permutations groups. Then we show that a 2-transitive locally
compact group (which satisfies some additional hypotheses) is
automatically a Lie group. Up to this point, our proof is more
or less the same as Tits' original proof; these results can also
be found in Salzmann \emph{et al.} 96.15.
After this point, we follow a different line than Tits.
First we consider the case where the space $X$ which $G$ acts on is
compact. This case is much easier than the noncompact case, since
one can use a convenient criterion due to 
Furstenberg which characterizes parabolics in Lie groups.
Using this result, it is not difficult to show that $G$ is
essentially a simple (noncompact) Lie group, and that
$X$ is a vertex set of the building $\Delta$ belonging to $G$.
If the real rank of the group $G$ is at least $2$,
the building has to be a projective space, and this leads
to a full classification in the compact case.

The noncompact case is more involved. Here we use some algebraic
topology to prove that $X$ is contractible if it is noncompact.
Once this is proved, it is not difficult to see that $X$ is
a real vector space, and that $G$ acts through affine-linear
transformations. The task is then to classify linear Lie groups
acting transitively on nonzero vectors in a real vector
space. Such a group acts transitively on the half-rays in the
vector space, and we can use the classification of compact Lie groups
acting transitively on spheres. In the last section we classify
those 2-transitive groups in our list which are Moufang sets.
This re-proves and
generalizes the classification by Kalscheuer, Tits, and Grundh\"ofer
of all sharply 2-transitive Lie groups.

\tableofcontents

\section{Preliminaries}

In this section
we collect a few general and simple
facts about 2-transitive groups.
An \emph{action} of a group $G$ on a nonempty set $X$ is a homomorphism
$G\too\mathrm{Sym}(X)$ of $G$ into the symmetric group of
$X$. For $x\in X$, we denote the stabilizer of $x$ by
$G_x=\{g\in G\mid g(x)=x\}$. For a subset $U\subseteq G$
we put $U\cdot x=\{g(x)\mid g\in U\}$.
If $G$ acts transitively on $X$ 
(i.e. if $X=G\cdot x$ for some $x\in X$)
and if $H\leq G$ is a subgroup, then
$H$ acts transitively on $X$ if and only if $G=G_xH$.
The \emph{kernel} of an action is the subgroup
$G_{[X]}=\bigcap\{G_x\mid x\in X\}$. The action is
\emph{effective} if $G_{[X]}=1$ (often this is called a
\emph{faithful} action). If the action is effective, $G$ can
be identified with its image in $\mathrm{Sym}(X)$.
A transitive action is \emph{regular} if
$G_x=1$ holds for all $x\in X$.

The action of $G$ on $X$ is called \emph{$k$-transitive}, 
for $|X|\geq k\geq 2$, if $G$ acts transitively on
the set of all $k$-tuples
$(x_1,\ldots,x_k)\in X^k$ with pairwise distinct entries.
Clearly,
a $(k+1)$-transitive group acts also $k$-transitively; in particular,
it acts transitively on $X$. If $|X|\geq 3$, then
$G$ acts 2-transitively on $X$ if and only if the stabilizer $G_x$ of
$x\in X$ acts transitively on $X\setminus\{x\}$ for every $x\in X$.
A transitive action (of $G$ on $X$)
is \emph{primitive} if the stabilizer $G_x$
is a maximal subgroup of $G$, i.e. if $G_x<H\leq G$ implies $H=G$.
\begin{Lem}
\label{NormalTrs}
If $G$ acts primitively on $X$, and if $N\unlhd G$ is a normal subgroup,
then either $N\leq G_x$ (and thus $N\leq G_{[X]}$),
or $N$ acts transitively on $X$.

\proof
See eg. Dixon-Mortimer \cite{DM} Theorem 1.6A(v).
\qed
\end{Lem}
The following criterion for primitivity is very convenient.
Let $R\subseteq X\times X$ be an equivalence relation,
i.e. $xRyRz$ implies $xRz$,
$\id_X\subseteq R$, and $\iota(R)=R$, where $\iota(x,y)=(y,x)$.
Let us call $\id_X=\{(x,x)\mid x\in X\}$ and $X\times X$ the
\emph{trivial equivalence relations}. An equivalence relation is
$G$-invariant if $g(R)=R$ holds for all $g\in G$
(with respect to the diagonal action of $G$ on $X\times X$).
\begin{Lem}
\label{NoBlocks}
A transitive action of $G$ on $X$ is primitive if and only if
$X\times X$ contains no nontrivial $G$-invariant equivalence
relations.
\proof
See eg. Jacobson \cite{Jac} Theorem 1.12, or
Dixon-Mortimer \cite{DM} Corollary 1.5A.
\qed
\end{Lem}
If $G$ is 2-transitive, then $G$ acts transitively on
$Y=X\times X\setminus\id_X$,
hence there are no nontrivial $G$-invariant equivalence relations.
Thus, a 2-transitive action is primitive.
\begin{Lem}
\label{GenLem}
Suppose that $G$ acts 2-transitively on $X$. Then the stabilizer
$G_x$ is a maximal subgroup of $G$.
\qed
\end{Lem}

\begin{Lem}
\label{AbGp}
Suppose that $G$ acts primitively and effectively on $X$.
If $1\neq A\leq G$ is an abelian normal subgroup,
then $G$ is a semidirect
product $G=G_x\ltimes A$. The group $A$ acts regularly on $X$, and if
we identify $A$ with $X$ via the evaluation map $ev_x:a\mapstoo a(x)$,
then each $g\in G_x$ acts by conjugation, 
\[
ga(x)=gag^{-1}g(x)=gag^{-1}(x).
\]
In particular, $G_x$ acts transitively on $A\setminus\{1\}$ if
$G$ acts 2-transitively on $X$.

\proof
By Lemma \ref{GenLem}, $A$ acts transitively on $X$. A transitive
and effective abelian group acts regularly. Thus $A_x=G_x\cap A=1$, and
$G=G_xA$ is a semidirect product.
\qed
\end{Lem}
In the next lemma, we identify $G$ with its image in $\mathrm{Sym}(X)$.
The result is used in Section 3.

\begin{Lem}
\label{CenLem}
Suppose that $G$ acts primitively and effectively on $X$.
If $G$ is not cyclic
of prime order, $G\not\cong\ZZ/p$, then
the centralizer of $G$ in the symmetric group is trivial,
$\mathrm{Cen}_{\mathrm{Sym}(X)}(G)=1$. In particular,
$G$ is centerless, and the normalizer
$\mathrm{Nor}_{\mathrm{Sym}(X)}(G)$ is contained in the
automorphism group of $G$.

\proof
Let $C=\mathrm{Cen}_{\mathrm{Sym}(X)}(G)$, and assume that $C\neq 1$.
Then $C$ is a nontrivial normal
subgroup in the primitive group $C G\leq\mathrm{Sym}(X)$,
hence $C$ acts transitively on $X$ by Lemma \ref{GenLem}.
Let $g\in G_x$. Then $c(x)=cg(x)=gc(x)$ holds
for all $c\in C$, so $g$ fixes the orbit $C\cdot x=X$
elementwise. It follows that $G$ acts regularly on $X$. Since the
action is primitive, $G$ has no proper subgroups, thus
$G\cong\ZZ/p$, for some prime $p$.
\qed
\end{Lem}

\section{Two-transitive Lie groups and locally compact groups}

We fix some topological
terminology. A homeomorphism is denoted by '$\cong$', and
a homotopy equivalence by '$\simeq$'. \emph{Unless stated otherwise,
all spaces are assumed to be Hausdorff.}
An \emph{$n$-manifold} is a second countable metrizable space
which is locally homeomorphic to $\RR^n$. In a topological group
we assume always that the singletons are closed, so the groups
themselves are regular topological spaces.
A \emph{Lie group} is a second countable topological group
which is at the same time a smooth manifold, such that
the group operations are smooth maps.
A \emph{transformation group} is a pair
$(G,X)$ consisting of a topological
group $G$, acting as a group
of homeomorphisms on a topological space $X\neq\emptyset$, such that
the map $G\times X\too X$ is continuous. If $G$ is a Lie group,
we call $(G,X)$ a \emph{Lie transformation group}.
The connected component of the identity in a topological group $G$
is denoted $G^\circ$.

A space is called \emph{$\sigma$-compact} if it is a countable
union of compact subsets.
With our conventions, Lie groups are $\sigma$-compact.
If $(G,X)$ is a transitive transformation group, and if
$x\in X$, then the map $G\too G/G_x$ is open, and thus the natural
map $G/G_x\too X$ is continuous. We need the following 
sharpening of this fact, cp.~Salzmann \emph{et al.} \cite{CPP} 96.8.

\begin{Prop}
\label{OpenAct}
Let $(G,X)$ be a transitive transformation group, and let
$x\in X$. If $G$ and $X$ are locally compact, and if $G$ is
$\sigma$-compact, then the natural map
$G/G_x\too X$ is a homeomorphism; in particular, the evaluation map
map $ev_x:G\too X$, $g\mapstoo g(x)$ is open.
\qed
\end{Prop}
The Approximation Theorem states that every locally compact
group can be approximated by Lie groups in the following sense,
cp. Salzmann \emph{et al.} \cite{CPP} 93.8. Note that every open subgroup
of $G$ contains the connected component $G^\circ$, and that
$G^\circ$ is a closed normal subgroup of $G$. However,
$G^\circ$ need not be open if $G$ is not a Lie group.
\begin{Num}\textbf{\hspace{-1ex}Approximation Theorem}
\label{ApprThm}
\em Let $G$ be a locally compact group.

(1)
There exists an open (and $\sigma$-compact)
subgroup $V\leq G$ such that $V/G^\circ$
is compact.

(2) Given  a subgroup $V$ as in (1), and given any neighborhood
$U$ of $1\in G$, there exists a compact normal subgroup $N\unlhd V$
of $V$ with $N\subseteq U$, such that $V/N$ is a Lie group.
\qed
\end{Num}
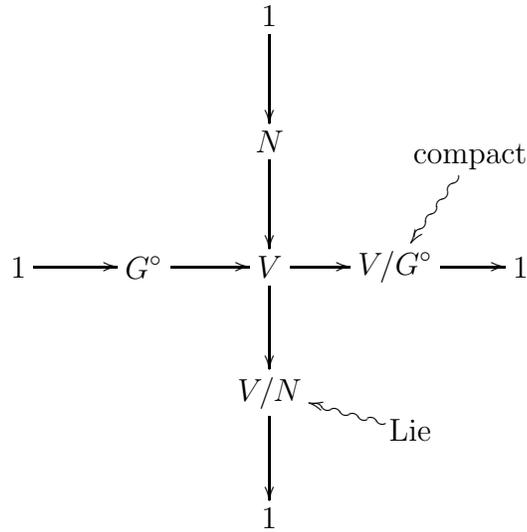
\[
\xymatrix@!=2em{
&& 1\ar[d] \\
&& N\ar[d] & {}
\save[]+<1cm,-.2cm>*\txt{compact}
\ar@{~>}[d]\restore \\
1\ar[r] & G^\circ\ar[r] & V \ar[r]\ar[d]
 & V/G^\circ\ar[r] & 1 \\
&& V/N\ar[d] &
\save[]+<.2cm,-.5cm>*\txt{Lie}
\ar@{~>}[l]\restore \\
&& 1}
\]
Theorem \ref{2-TrsLie} below is essentially due to Tits \cite{Ti};
the present proof is a simplified version of the proof
given in Salzmann \emph{et al.} 
\cite{CPP} 96.15.
A space is called \emph{totally disconnected} if the component of
every point is trivial.
The following useful result is due to Eilenberg \cite{Eil} 3.1.
\begin{Lem}
\label{EilLem}
Let $X$ be a connected space, and let
$Y=X\times X\setminus\id_X=\{(x,y)\mid x,y\in X,\,x\neq y\}$.
Let $\iota(x,y)=(y,x)$.
If $Y$ is not connected, then $Y$ has precisely two 
(necessarily open) components $A,B$, and $\iota(A)=B$.
\qed
\end{Lem}

\begin{Thm}
\label{2-TrsLie}
Let $(G,X)$ be an effective topological transformation group,
where $G$ is $\sigma$-compact and locally compact, and $X$ is
locally compact and not totally disconnected. Suppose
that $G$ acts 2-transitively on $X$.
Then $G$ is a Lie group and $X\cong G/G_x$ is a smooth and
connected manifold. Every open subgroup $V\leq G$ acts primitively on
$X$; in particular, $G^\circ$ acts primitively on $X$.
The group $G^\circ$ is noncompact and has at most three orbits
in $X\times X$.

\proof
First we prove that $X$ is connected. Define an equivalence relation
$R$ on $X$ by setting $xRy$ if $x,y$ are contained in a connected subset.
Since $X$ is not totally disconnected, $R\neq\id_X$. This relation is
$G$-invariant, since $G$ acts by homeomorphisms on $X$.
Thus $R=X\times X$ by Lemma \ref{NoBlocks}, and hence $X$ is connected.

Now we show that $G^\circ\neq 1$. If $G^\circ=1$, then $G$ is
totally disconnected, and thus zero-dimensional, see
Hewitt-Ross \cite{HR} Thm.~3.5. Therefore, $X\cong G/G_x$ is also
zero-dimensional, see \emph{loc.cit.} Thm.~7.11, and $X$
contains arbitrarily small open and closed subsets. This contradicts
the fact that $X$ is connected.

Let $Y=\{(x,y)\mid x,y\in X,\,x\neq y\}$. This is a locally
compact space, and $G$ acts transitively on $Y$; the
evaluation map $ev_{(x,y)}:G\too Y$ is open by Proposition
\ref{OpenAct}. Let $V\leq G$ be an open subgroup.
Since $V$ contains the normal subgroup $G^\circ$, it acts
transitively on $X$ by Lemma \ref{NormalTrs}.
For every $(x,y)\in Y$, the orbit $ev_{(x,y)}(V)=V\cdot(x,y)$ is open.
By Lemma \ref{EilLem}, $Y$ has at most two components $A,B$, so
$V$ has at most two orbits $A,B$ in $Y$. Moreover, $\iota(A)=B$,
so $X\times X$ contains no nontrivial $V$-invariant
equivalence relation, i.e. the action of $V$ on $X$
is primitive by Lemma \ref{NoBlocks}.

Now choose a cocompact open subgroup $V\leq G$ 
as in the Approximation Theorem
\ref{ApprThm} (1), and choose a small neighborhood $U\subseteq V$
of the identity such
that $U\cdot x\neq X$. Then $U$ cannot contain a proper normal
subgroup $N\unlhd V$, since otherwise $N\cdot x=X$.
Therefore, $V^\circ=G^\circ$ is an open (Lie) subgroup of
$G$ by Theorem \ref{ApprThm},
and $G/G^\circ$ is discrete and (by $\sigma$-compactness
of $G$) countable.

If the open subgroup $G^\circ$ is compact, then its orbits in
$X\times X$ are compact. Since there are at most three different
orbits $A,B,\id_X$, this would imply
that $A\cup B=X\times X\setminus\id_X$ is compact, a contradiction to the
fact that $X\setminus\{x\}$ is noncompact (because $X$ is connected).
\qed
\end{Thm}
\emph{For the remainder of this section,
we assume that $G$ is a Lie transformation group
acting on a connected manifold $X\cong G/G_x$,
and that this action is effective and 2-transitive.}

We need two results about the connectivity of complements of discrete
subsets.

\begin{Lem}
Let $M$ be a connected manifold of positive dimension.
If $M$ is compact or if $\dim(M)\geq 2$, then 
$M\setminus\{x\}$ is path connected for all $x\in M$.

\proof
If $\dim(M)\geq 2$, then $H_1(M,M\setminus\{x\})=0=\widetilde H_0(M)$,
and the exact sequence 
\[
\too H_1(M,M\setminus\{x\})\too
\widetilde H_0(M\setminus\{x\})\too \widetilde H_0(M)\too 0
\]
shows that
$\widetilde H_0(M\setminus\{x\})=0$, so $M\setminus\{x\}$ is also
path connected.
If $\dim(M)=1$ and if $M$ is compact, then $M\cong\SS^1$ is a circle,
and $M\setminus\{x\}\cong\RR$ is path connected.
\qed
\end{Lem}

\begin{Lem}
\label{CC2trs}
If $X$ is compact or if $\dim(X)>1$, then
the connected component $G^\circ$ acts also 2-transitively
on $X$.

\proof
The connected component $G^\circ\leq G$ is a normal
subgroup, hence $G^\circ$ acts transitively on $X$ by 
Lemma \ref{GenLem}.
The subgroup $G^\circ\leq G$ is open, hence
$(G^\circ)_x$ is open in $G_x$. Let $y\in X\setminus\{x\}$.
The evaluation map 
$ev_y:g\mapstoo g(y)$, $G_x\too X\setminus\{x\}$ is open by
Proposition \ref{OpenAct};
since $X\setminus\{x\}$ is connected,
$ev_y((G^\circ)_x)=X\setminus\{x\}$, i.e. $(G^\circ)_x$ acts transitively
on $X\setminus\{x\}$. 
\qed
\end{Lem}

\begin{Lem}
\label{AbLem}
If $G$ has an abelian normal subgroup $A\neq 1$, then $A\cong\RR^n$ is a
vector group.

\proof
The group $A$ acts regularly on $X$ by Lemma \ref{AbGp},
hence $A$ is connected (by Proposition \ref{OpenAct}).
A connected abelian $n$-dimensional Lie group is of the form
$\RR^k/\ZZ^k\times\RR^{n-k}$, for some $k\in\{0,\ldots,n\}$,
see Onishchik-Vinberg \cite{OnVin} Ch.~1 \S2 Prop.~3, p.~29.
If $k\geq 1$, then $\RR^k/\ZZ^k$
contains elements of order $l$ for all $l\in\mathbb{N}$.
This is not possible by Lemma \ref{AbGp} (all elements of $A$
different from the identity
have to be conjugate under $G_x$), so $k=0$.
\end{Lem}
In particular, the existence of a nontrivial proper abelian normal
subgroup implies that $X$ is noncompact.
The converse will be proved in Section 4: if $X$ is noncompact, then
$G$ has a nontrivial proper abelian normal subgroup.

\begin{Prop}
\label{SimpleLem}
If $G$ has no nontrivial proper normal abelian subgroup, then
the connected component $G^\circ$ is simple. 

\proof
Assume that $G$ has no nontrivial proper abelian normal subgroup.
Let $\sqrt G$ denote the radical of $G$, i.e. the largest normal
connected solvable subgroup, see Onishchik-Vinberg
\cite{OnVin} Ch.~1 \S4 6$^\circ$. Let $1<Z<\cdots<\sqrt G$ denote the
derived series of $\sqrt G$. If $\sqrt G\neq 1$, then $Z$ is a
nontrivial normal abelian
subgroup in $G$, since it is characteristic in $\sqrt G$. Thus $\sqrt G=1$.
Note also that $\sqrt G=\sqrt{G^\circ}$. Therefore $G^\circ$ is semisimple
and centerless and thus a direct product of simple Lie groups,
see Salzmann \emph{et al.} \cite{CPP} 94.23.

If $\dim(X)\geq 2$ or if $X$ is compact, then $G^\circ$ acts
2-transitively on $X$ by Lemma \ref{CC2trs}.
If $\dim(X)=1$ and if $X$ is noncompact, then
$X\cong\RR$. By Brouwer's result, no centerless semisimple Lie group
acts transitively on $\RR$, see Salzmann \emph{et al.} \cite{CPP} 96.30,
so this case cannot occur.

Assume now that $\dim(G)\geq 2$, and that
$G^\circ=G_1\times G_2$, with $G_1\neq 1$
simple. Suppose that $G_2\neq 1$.
Then $G_1$ and $G_2$ are normal in $G^\circ$ and  thus
transitive on $X$. Suppose that
$g_1\in G_1$ fixes $x$, and that $g_2\in G_2$.
Then $g_2(x)=g_2g_1(x)=g_1(g_2(x))$, so $g_1$ fixes $g_2(x)$
for all $g_2\in G_2$. Since $G_2\cdot x=X$, this implies
that $g_1=1$, i.e. $G_1$ acts regularly on $X$. We identify $X$ with
$G_1$; then the action of $G_1$ on $X$ is the standard left regular
action $(g,x)\mapstoo gx$ of $G_1$ on itself. The centralizer of this action
$\mathrm{Cen}_{\mathrm{Sym}(X)}(G_1)$ is isomorphic to $G_1$
with the action $(g,x)\mapstoo xg^{-1}$. 
It follows that $G_2\cong G_1$, with the action
$(g_1,g_2)(x)=g_1xg_2^{-1}$. The stabilizer of $1\in X$ is the
diagonal subgroup $\{(g,g)\mid g\in G_1\}$; its orbits in $X$ are
the conjugacy classes of $G_1$. But a centerless simple Lie group
contains a torus $\RR/\ZZ$; in particular, it contains elements
of arbitrary finite order, and thus infinitely many conjugacy classes.
This is a contradiction to the 2-transitivity of $G^\circ$.
Therefore $G_2=1$ and $G=G_1$ is simple.
\qed
\end{Prop}
If we assume in addition that $X$ is compact, then there is a
shortcut in the last step of the proof above:
if $G^\circ$ is a product of at least two simple factors, then each
simple factor acts regularly, and thus $G^\circ$ is compact,
a contradiction to Theorem \ref{2-TrsLie}.

\section{The case when $X$ is compact}
\label{CompactCase}
In this section we classify all pairs $(G,X)$, where $G$ is a
Lie group acting effectively and
2-transitively on a compact connected manifold $X$.
By \ref{CC2trs}, \ref{AbLem}, and \ref{SimpleLem},
$G^\circ$ is a simple centerless Lie group which
acts 2-transitively on $X$.
Recall the definition of a parabolic subgroup in an arbitrary
semisimple Lie group $H$ (of noncompact type). Suppose that 
\[
H=KAU
\]
is an Iwasawa decomposition for $H$ (the unipotent subgroup
$U$ is denoted '$N$' in most books on Lie groups; we use $U$, since
$N$ has a different meaning with BN-pairs), see Helgason
\cite{He} Ch.~IX  \S1 Thm.~1.3.
Let $K_0=\mathrm{Cen}_K(A)$ denote the $K$-centralizer of $A$
(this is the \emph{reductive anisotropic kernel} of $H$, see
Tits \cite{TitsAG} p.~39). Then
\[
B=K_0AU
\]
is a subgroup, a \emph{minimal parabolic subgroup} of $H$, see Helgason
\cite{He} Ch.~IX \S1, Warner \cite{War} p.~56.
A \emph{parabolic subgroup} is any subgroup of $H$ which is conjugate
to an overgroup of $B$. There are precisely
$2^{\dim(A)}$ conjugacy classes
of parabolics in $H$ (including the classes of $H$ and $B$), see Warner
\cite{War} Thm.~1.2.1.1. If we put $N=\mathrm{Nor}_H(K_0A)$, then
$(B,N)$ is a BN-pair of rank $\dim(A)$ for the group $H$, see
Warner \cite{War} p.~68. 

\begin{Thm}
\label{FuThm}
Let $H\leq\mathrm{PSL}_{m+1}\RR$ be a closed semisimple subgroup, 
and consider the standard projective action of $\mathrm{PSL}_{m+1}\RR$
on the real projective space
$\RR\mathrm{P}^m$. Let $X\subseteq\RR\mathrm{P}^m$ be a closed
$H$-invariant subset.
Assume that the following two conditions are satisfied.

(a) If $\emptyset\neq Y\subseteq X$ is closed and $H$-invariant, then
$Y=X$.

(b) If $\emptyset\neq Z\subseteq X\times X$ is closed and
$H$-invariant, then $Z\cap\id_X\neq\emptyset$.

\noindent
Then $X=H/H_x$ is a homogeneous space of $H$, and $H_x$ is a
parabolic subgroup of $H$.

\proof
Condition (a) says that the action is \emph{minimal} in the sense
of Furstenberg \cite{Fu} p.~278--279, and condition (b) says that it
is \emph{proximal}; by the very definition, it is what Furstenberg
calls \emph{projective}. By Proposition 4.3. in
\emph{loc.cit.} p.~280, there exists an $H$-equivariant map
from the \emph{maximal Furstenberg boundary} $B(H)$ of $H$ onto $X$.
Since we assumed that $H$ is semisimple, the discussion
in \emph{loc.cit.} p.~280 shows that there is an $H$-equivariant
homeomorphism $B(H)\cong H/B$ between the Furstenberg boundary and
the coset space $H/B$, where $B\leq H$ is a minimal parabolic. The space
$H/B\cong K/K_0$ is compact (it is the chamber set of the canonical
building $\Delta(H)$ associated to $H$), and thus
$H_x$ contains a conjugate of $B$,  i.e. $H_x$ is a parabolic.
\qed
\end{Thm}

\begin{Prop}
\label{StabIsPara}
If $X$ is compact, and if $G$ is a simple connected Lie group acting
effectively and 2-transitively on $X$, then
the stabilizer $P=G_x$ is a maximal parabolic subgroup
of the simple group $G$.

\proof Let $\fp$ denote the Lie algebra of $P$, with
$k=\dim_\RR(\fp)$.
The $G$-normalizer of $\fp$ contains $P$; by maximality
of $P$, it coincides with $P$. Consider the adjoint
action of $G$ on the Grassmann manifold $\mathrm{Gr}_k(\fg)$
of all $k$-dimensional real subspaces of the Lie algebra $\fg$
of $G$.
The stabilizer of $\fp\in\mathrm{Gr}_k(\fg)$ is $P$, as we just saw. 
The Grassmannian embeds $G$-equivariantly into the real projective
space $\mathbf{P}(\bigwedge^k\fg)$ of the $k$-th exterior
power of the real vector space $\fg$. The $G$-orbit
of $\bigwedge^k\fp\in\mathbf{P}(\bigwedge^k\fg)$ 
is therefore $G$-equivariantly homeomorphic to $X=G/P$.
Since $G$ acts 2-transitively, the conditions (a) and (b) in
Furstenberg's Theorem \ref{FuThm} are clearly satisfied, hence
$P$ is a parabolic subgroup, and, by 2-transitivity on $X$, a
maximal one.
\qed
\end{Prop}
Recall that the \emph{real rank} of a semisimple
Lie group $H$ is the real vector space dimension of $\mathfrak a$,
where $\mathfrak{h}=\mathfrak{k}+\mathfrak{a}+\mathfrak{u}$
is an Iwasawa decomposition of the Lie algebra $\mathfrak{h}$ of
$H$.

\begin{Thm}
\label{CompactClass}
Let $(G,X)$ be an effective 2-transitive transformation group.
Assume that $G$ is locally compact and $\sigma$-compact, and
that $X$ is compact and not totally disconnected. Then $G$ is a
Lie group and the connected component $G^\circ$ is simple and
noncompact. There exists a subgroup $\widehat G$ of the
automorphism group of $G^\circ$ such that either $G=G^\circ$,
or $G=\widehat G$. Both $\widehat G$ and $G^\circ$ act 2-transitively
on $X$. Topologically, the space $X$ is either a sphere or a projective
space over $\RR$, $\CC$, $\HH$, or $\OO$.
The possibilities for the pair $(G^\circ,X)$ are as follows.

\medskip
(a) If the real rank of $G^\circ$ is $1$, then $X$ is a sphere,
and there are the following possibilities.

\begin{center}
\begin{tabular}{c|c|c|cl} 
$G^\circ$ & $X$ & $\widehat G$ & $|\widehat G/G^\circ|$ \\ \hline
$\mathrm{PEO}_{n,1}\RR$  & $\SS^{n-1}$  & 
$\mathrm{PO}_{n,1}\RR$ & $2$ & ($n\geq 2$) \\
$\mathrm{PSU}_{n,1}\CC$ & $\SS^{2n-1}$ & 
$\mathrm{P}\Gamma\U_{n,1}\CC$ & $2$ & ($n\geq 2$) \\
$\mathrm{PU}_{n,1}\HH$ & $\SS^{4n-1}$ & 
$\mathrm{PU}_{n,1}\HH$ & $1$ & ($n\geq 2)$\\
$\mathrm{F}_{4(-20)}$        & $\SS^{15}$   & 
$\mathrm{F}_{4(-20)}$ & $1$ 
\end{tabular}
\end{center}

\medskip
\noindent
The simple subgroup
$\mathrm{PEO}_{n,1}\RR=(\mathrm{PO}_{n,1}\RR)^\circ$
is the subgroup generated by all Eichler (or Siegel) transformations,
see e.g. Hahn-O'Meara \cite{HOM} p.~214. Alternatively, the group
$\mathrm{PEO}_{n,1}\RR$ can be described as the commutator group of
$\mathrm{PO}_{n,1}\RR$, or as the connected component
$(\mathrm{PO}_{n,1}\RR)^\circ$.

\medskip
(b) If the real rank $k$ of $G^\circ$ is at least $2$, then 
$X$ is a projective space of rank $k$, and there are only
the following possibilities.

\begin{center}
\begin{tabular}{c|c|c|cl} 
$G^\circ$ & $X$ & $\widehat G$ & $|\widehat G/G^\circ|$ \\ \hline
$\mathrm{PSL}_{k+1}\RR$ & $\RR\mathrm{P}^k$ & 
$\mathrm{PGL}_{k+1}\RR$ & $2$ & ($k$ odd)\\
$\mathrm{PSL}_{k+1}\RR$ & $\RR\mathrm{P}^k$ & 
$\mathrm{PSL}_{k+1}\RR$ & $1$ & ($k$ even)\\
$\mathrm{PSL}_{k+1}\CC$ & $\CC\mathrm{P}^k$ & 
$\mathrm{P}\Gamma\mathrm{L}_{k+1}\CC$ & $2$ \\
$\mathrm{PSL}_{k+1}\HH$ & $\HH\mathrm{P}^k$ & 
$\mathrm{PSL}_{k+1}\HH$ & $1$ \\
$\mathrm{E}_{6(-26)}$       & $\OO\mathrm{P}^2$ & 
$\mathrm{E}_{6(-26)}$ & $1$ & ($k=2$) 
\end{tabular}
\end{center}
\end{Thm}
Note that there are the following isomorphisms of simple
Lie groups:
\begin{gather*}
\mathrm{PSL}_2\RR\cong\mathrm{PSU}_{1,1}\RR\cong\mathrm{PEO}_{2,1}\RR\\
\mathrm{PSL}_2\CC\cong\mathrm{PEO}_{3,1}\RR\\
\mathrm{PU}_{1,1}\HH\cong\mathrm{PEO}_{4,1}\RR\\
\mathrm{PSL}_2\HH\cong\mathrm{PEO}_{5,1}\RR
\end{gather*}

\medskip
\emph{Proof of Theorem \ref{CompactClass}.}
We know already by \ref{CC2trs}, \ref{AbLem}, and \ref{SimpleLem}
that $H=G^\circ$ is a simple centerless
Lie group which acts 2-transitively on $X$, and by Proposition
\ref{StabIsPara} the stabilizer $P=H_x$ is a maximal parabolic
subgroup. If the real rank of $H$ is 1, then
Table V in Ch.~X, p.~518 in Helgason \cite{He}
shows that $(H,X)$ is a as in (a). Note that a group with a BN-pair
of rank 1 is the same as a 2-transitive group; thus, each simple
Lie group of real rank 1 is a 2-transitive group. Of course, this
follows also by direct inspection of the list of groups.

Now suppose that the real rank $k$ of $H$  is at least 2.
Then we consider the irreducible
spherical building $\Delta=\Delta(H)$ of $H$. Theorem \ref{2transbuild}
in the next
section shows that $\Delta$ is the building associated to a projective
space (of finite rank, as all buildings related to Lie groups are
of finite rank),
and that the parabolic $H_x$ is the stabilizer of a point
in this projective space.

Then $H$ is one of the
groups $\mathrm{PSL}_{k+1}\FF$, for $\FF=\RR,\CC,\HH$, and
$\Delta$ is a projective geometry of rank $k$ over $\FF$, or
$k=2$, and $\Delta$ is the Cayley plane and $H=\mathrm{E}_{6(-26)}$.
This follows either from the topological Fundamental Theorem of
Projective Geometry, (as proved by Kolmogorov \cite{Kolm}, see
also K\"uhne-L\"owen \cite{KueLoe} and the survey by
Grundh\"ofer-L\"owen \cite{GruLoe}),
or from Cartan's classification of simple Lie groups
and Tits' classification \cite{TitsAG} of their Tits diagrams 
(unjustly sometimes called ''Satake diagrams'')
and their relative diagrams, see
Helgason \cite{He} Table VI, Ch.~X p.~532. In any case,
we obtain the list of groups $G^\circ$ and spaces $X$ in (b).

By Lemma \ref{CenLem}, $G$ is a subgroup of the automorphism group of
$H=G^\circ$. The outer automorphism groups of all simple
Lie groups have been determined by Takeuchi \cite{Ta},
see also Wolf \cite{Wolf} 8.8 p.~263.
In case (b), we have the additional condition that $G$ preserves
the conjugacy class of $P$, and thus $G$ has to act by type-preserving
automorphisms on the projective space. Thus, we obtain the list
of groups $\widehat G$ given in (a) and (b).
\qed

\begin{Cor}
If $G$ is 3-transitive, then $X\cong\SS^{n-1}$
and $G\cong\mathrm{PO}_{n,1}\RR$, for $n\geq 2$, or
$G\cong\mathrm{PEO}_{n,1}\RR$ and $n\geq 3$.
None of these groups is 4-transitive.

If the action is sharply 3-transitive (i.e. regular on triples
of pairwise distinct points), then $G\cong\mathrm{PGL}_2\RR$ or
$G\cong\mathrm{PSL}_2\CC$.

\proof
If the building $\Delta=\Delta(G^\circ)$
is a projective geometry of rank $k\geq 2$,
then there exists triples of collinear points and triples of
noncollinear points, so $G$ cannot be 3-transitive. Therefore,
a 3-transitive group is one of the groups given in 
Theorem \ref{CompactClass} (a).
The unitary groups $\mathrm{P}\Gamma\U_{n,1}\CC$ 
and $\mathrm{PU}_{n,1}\HH$ cannot act 3-transitively if $n\geq 2$:
choose three linearly independent vectors $u,v,w$ with
$(u|v)=(v|w)=(w|u)$ and $(u|u)=(v|v)=(w|w)=0$
(we denote the corresponding hermitian form by
$(x|y)=-\bar x_0y_0+\bar x_1y_1+\cdots+\bar x_ny_n$). Then
$x=u+v\ti$ is also isotropic, but no semilinear map fixing
the subspace spanned by $\{u,v\}$ can move the subspace spanned
by $x$ to the subspace spanned by $w$. There is a more geometric
way to see this, which is also valid for
the group $\mathrm{F}_{4(-20)}$. These groups can be viewed
as subgroups of $\mathrm{P}\Gamma\mathrm{L}_{n+1}\CC$,
$\mathrm{PGL}_{n+1}\HH$, acting on the set $\mathcal Q$ of absolute
points of a hyperbolic polarity in the projective geometry
(for the group $\mathrm{F}_{4(-20)}$, the corresponding projective
space is the Cayley plane $PG_2(\mathbb{O})$, see \cite{CPP} Ch.~1).
If we fix two distinct points $x,y\in\mathcal Q$, then
we fix the line $L$ joining these two points. Now $L\cap\mathcal O$
is a proper subset of $\mathcal Q$, and
$L\cap{\mathcal Q}\neq \{x,y\}$ if $n\geq 2$; in fact,
$L\cap{\mathcal Q}\cong \SS^{\dim_\RR\FF-1}$ for
$\FF=\RR,\CC,\HH,\mathbb{O}$. This excludes the unitary groups
and $\mathrm{F}_{4(-20)}$. 

The group
$\mathrm{PEO}_{2,1}\RR\cong\mathrm{PSL}_2\RR$ is not 3-transitive;
the remaining orthogonal groups are 3-transitive, as can be checked.
Finally, none of these groups is 4-transitive: Let $p,q,r,s\in\mathcal Q$
be four distinct points such that the lines $pq$ and $rs$ intersect
in an interior point of $\mathcal Q$. Then $r,s$ cannot be moved to two
points $r',s'$ (fixing $p,q$) such that $r's'$ intersects $pq$ in an
exterior point. If the action is sharply 3-transitive, then
$\dim(G)=3\dim(X)$, and thus only the groups $\mathrm{PGL}_2\RR$ and
$\mathrm{PSL}_2\CC$ remain.
\qed
\end{Cor}

\begin{Cor}
If $X$ is a complex homogeneous space (i.e. if $X$ is a complex
manifold and if $G$ preserves the complex structure), then
$X=\CC\mathrm{P}^n$ and $G\cong\mathrm{PSL}_{n+1}\CC$..

\proof
Suppose that $X$ is a complex manifold, and that the $G$-action
preserves the complex structure. Then $\dim(X)$ is even.
In case (a) of Theorem \ref{CompactClass},
$X\cong\SS^{2k}$ is a sphere, and
$G^\circ=\mathrm{PEO}_{2k+1,1}\RR$. The compact subgroup
$K=\SO(2k+1)$ acts transitively on $X$, and the
isotropy representation of $K_x$ on the tangent space
$T_xX\cong\RR^{2k}$ is the standard action of $\SO(2k)$.
If $k\geq 2$, then this action is not complex linear.
Thus $k=1$ and $G^\circ=\mathrm{PEO}_{3,1}\RR\cong\mathrm{PSL}_2\CC$.

In case (b) of \ref{CompactClass},
$X$ is an even-dimensional real projective space or
a complex or quaternionic projective space. Again we consider
the isotropy representation of $K_x$ for a suitable compact subgroup
$K\leq G^\circ$. For $\RR\mathrm{P}^{2k}$, we have $K=\SO(2k+1)$,
and $K_x=\mathrm{O}(2k)$ acts in the standard way on $T_xX\cong\RR^{2k}$.
This action preserves no complex structure on the tangent space, not
even for $k=1$.
For $X=\HH\mathrm{P}^n$ we take $K=\Sp(n+1)$; then
$K_x=\Sp(n)\cdot\Sp(1)$ acts in the standard way
(from the left and right) on $T_xX\cong\HH^n$; again, there is
no invariant complex structure for all $k\geq 1$.
Thus the only possibility is that $G^\circ=\mathrm{PSL}_{n+1}\CC$
in the standard action on $\CC\mathrm{P}^n$. Complex conjugation
is not $\CC$-linear, hence $G=G^\circ$.
\qed
\end{Cor}

\section{Digression: two-transitive actions on buildings}

In this section we consider the following situation. $\Delta$ is a
building (thick, and
of arbitrary, possibly infinite rank), and $G$ is a group of
(type-preserving) automorphisms of $\Delta$. We assume that the action of
$G$ on one type of residues of $\Delta$ is 2-transitive.
This is a purely combinatorial problem, and we make no topological
assumptions. After I had finished a first version of this section,
Bernhard M\"uhlherr pointed out to me that the book by
Brouwer-Cohen-Neumaier \cite{BCN} contains related results.

Recall that a \emph{Coxeter diagram} is
an undirected labeled
graph whose nodes are labeled by some set $I$. For each
pair $i,j\in I$, there is a number $m_{ij}=m_{ji}\in\mathbb{N}\cup\{\infty\}$,
with $m_{ii}=1$, and $m_{ij}\geq 2$ if $i\neq j$.
If the nodes $i,j$ are not adjacent (joined by an edge), then $m_{ij}=2$.
If the nodes $i,j$ are adjacent, and if $m_{ij}\geq 4$, then the
edge joining them is labeled with the number $m_{ij}$. 
The corresponding \emph{Coxeter group} is the group generated by a set
$S=\{s_i\mid i\in I\}$ of generators labeled by $I$,
with relations $(s_is_j)^{m_{ij}}=1$ (for
$m_{ij}\neq\infty$). The pair $(W,S)$ is called
a \emph{Coxeter system}. The Coxeter system is \emph{irreducible} if 
its Coxeter graph is connected. Reducible Coxeter systems are products
(in a natural sense). For $J\subseteq I$, the
subgroup generated by $S_J=\{s_J\mid j\in J\}$ is denoted $W_J$;
one can prove that the pair $(W,S_J)$ is again a Coxeter system.
We refer to the books by Bourbaki \cite{Bourb}
and Humphreys \cite{Humph}, and to
Ch.~2 in Ronan \cite{Ronan}.
Recall the Coxeter diagram $\mathbf A_k$ which is defined as
\begin{diagram}[abut,size=2em]
\underset{\rlap{$\scriptstyle 1$}}\bullet & \rLine &
\underset{\rlap{$\scriptstyle 2$}}\bullet & \rLine &
\underset{\rlap{$\scriptstyle 3$}}\bullet & \rLine &
{} & \cdots & {} & \rLine &
\underset{\rlap{$\scriptstyle k-2$}}\bullet & \rLine &
\underset{\rlap{$\scriptstyle k-1$}}\bullet & \rLine &
\underset{\rlap{$\scriptstyle k$}}\bullet
\end{diagram}
We denote the limit of these Coxeter diagrams (as $k$ goes to infinity)
by $\mathbf A_\omega$; the corresponding Coxeter diagram is thus 
\begin{diagram}[abut,size=2em]
\underset{\rlap{$\scriptstyle 1$}}\bullet & \rLine &
\underset{\rlap{$\scriptstyle 2$}}\bullet & \rLine &
\underset{\rlap{$\scriptstyle 3$}}\bullet & \rLine &
\underset{\rlap{$\scriptstyle 4$}}\bullet & \rLine &
\underset{\rlap{$\scriptstyle 5$}}\bullet & \rLine &
\underset{\rlap{$\scriptstyle 6$}}\bullet & \rLine & {} & \cdots & {}
\end{diagram}
The next result can also be extracted from Cooperstein \cite{Coop},
cp.~Brouwer-Cohen-Neumaier \cite{BCN} Thm.~10.2.3 p.~300.
\begin{Prop}
\label{DoubleCox}
Let $(W,\{s_i\mid i\in I\})$ be an irreducible Coxeter system, where
$I$ has (possibly infinite) cardinality $\kappa$. Let $J$ be
a subset of $I$ and assume that $W$ acts 2-transitively on $W/W_J$.
Then $\kappa\leq\aleph_0$, and either
the Coxeter system is of type $\mathbf A_k$, with $J=\{1\}$ or
$J=\{k\}$, or of type $\mathbf A_\omega$, and
$J=I\setminus\{1\}$
(the nodes of the diagram are labeled as above).

\proof
The proof is very simple.
We divide it into four basic steps.

(1) \emph{$J$ has corank $1$, i.e. $I\setminus J$ is a singleton.}

Since $W_J$ has to be a maximal subgroup, $I\setminus J$ is a singleton,
which we denote by $\{1\}$. Thus we have $W=W_J\cup W_J s_1 W_J$.
Now we use the standard description of shortest double coset
representatives: every element $w$ which has the property that
all reduced expressions of this
element start and end with $s_1$ represents a unique double coset
$W_J w W_J$.

(2) \emph{In the Coxeter diagram, the node $1$ has at most one neighbor.}

Assume that node $1$ has two different neighboring nodes, say $2,3$.
Then $s_1s_2s_3s_1$ is a reduced word representing a third double coset,
$W_Js_1s_2s_3s_1W_J$ (note that the order of $s_2s_3$ is not important
for this argument).

(3) \emph{No node in the Coxeter digram has more than two neighbors.}

Assume otherwise; let $i\neq 1$ be a node with three different neighbors.
Let $1$-$2$-$3$-$\cdots$-$(i-1)$-$i$
be a shortest path in the Coxeter diagram, and
let $j,k$ be distinct neighbors of node $i$ different from $i-1$. Then
$s_1s_2\cdots s_{i-1}s_is_js_ks_is_{i-1}\cdots s_1$ is a reduced word
which represents another shortest double coset representative.
This implies already that $I$ is finite or countable, and that the
Coxeter diagram of $W$ is a string diagram (with at least one end).

(4) \emph{The Coxeter diagram has no multiple bonds, $m_{ij}\leq 3$.}

Assume that $1$-$2$-$\cdots$-$i$-$j$
is a shortest path in the diagram, and that
$i,j$ are joined by an edge labeled $m_{ij}\geq 4$. Then $s_is_j$ has order at
least $4$. The reduced word
$s_1s_2\cdots s_{i-1}s_is_js_is_{i-1}\cdots s_1$ represents another shortest
double coset representative.

By (1) -- (4),
the diagram has only simple bonds and does not branch, and $1$ is an
end node.
\qed
\end{Prop}
Now we consider 2-transitive actions on buildings. Since we
allow infinite rank, we view buildings
as \emph{chamber systems}, see Tits \cite{TitsGV} or
Ronan \cite{Ronan}.
A chamber system is a set $\mathcal C$, endowed with a
collection $\{{\sim_i}\mid i\in I\}$ of equivalence relations.
For a subset $J\subseteq I$, the equivalence relation generated
by $\{{\sim_j}\mid j\in J\}$ is denoted $\sim_J$.
Given a Coxeter system $(W,\{s_i\mid i\in I\})$,
a \emph{building}
\[
\Delta=(\mathcal C,\{\sim_i\mid i\in I\},\delta,W,\{s_i\mid i\in I\})
\]
is a chamber system $\mathcal C$, endowed with a 
metric $\delta:\mathcal C\times \mathcal C\too W$ taking its
values in a Coxeter group (subject to
certain axioms which can be found eg. in Ronan's book \cite{Ronan}
Ch.~3).
An automorphism of a building (in this sense) is a permutation
of $\mathcal C$ which preserves the equivalence classes (and
thus the metric $\delta$). A $J$-\emph{residue} in a building
is a $J$-equivalence class.

For $i\in I$, let $V_i$ denote the set of all residues of type
$I\setminus\{i\}$ in the building. Define an \emph{incidence relation $*$}
as follows: two residues are incident if their intersection is
nonempty. The datum $\Gamma=(\{V_i\}_{i\in I},*)$ determines a 
\emph{diagram geometry}. Suppose now that the building is of type
$\mathbf A_\kappa$, for $2\leq\kappa\leq\omega$.
Put $\cP=V_1$ and $\cL=V_2$ (the nodes are labeled as before).
\begin{Thm}
If $\kappa\geq 2$, then the geometry $(\cP,\cL,*)$ obtained from
an $\mathbf A_\kappa$-building is a projective space.

\proof
The proof indicated by Tits in \cite{TitsGV} p.~540 applies (also for
the case $\kappa=\omega$).
\qed
\end{Thm}

\begin{Thm}
\label{2transbuild}
Let $\Delta$ be a thick irreducible building over some type set $I$,
let $G$ be a group of type preserving
automorphisms of $\Delta$ acting 2-transitively on the collection
of all $J$-residues, for some $J$.
Then $\Delta$ is of type $\mathbf A_k$ or $\mathbf A_\omega$;
the corresponding diagram geometry can be identified with a projective space
(possibly of infinite rank), and $G$ acts 2-transitively on the points
of this space.

\proof
Let $X,Y,Z$ be distinct $J$-residues. Since the stabilizer of $X$
is a maximal subgroup, $J$ has corank $1$, i.e. $I\setminus J$ is
a singleton. There is a unique element
$w\in\delta(X\times Y)\subseteq W$
such that $\ell(w)$ minimizes the lengths $\ell(\delta(x,y))$,
where $x,y$ run through all chambers in $X,Y$, see Scharlau-Dress \cite{DS}.
Since the action is 2-transitive
on the $J$-residues, the same $w$ gives the minimal distance for
elements in $Z$. Thus $W=W_J\cup W_JwW_J$; in particular, $W$ acts
2-transitively on $W/W_J$. By Proposition \ref{DoubleCox},
$\Delta$ is of type $\mathbf A_k$ or $\mathbf A_\omega$.
\qed
\end{Thm}

\section{The case when $X$ is noncompact}

In this section we assume that $G$ is a 2-transitive effective Lie
transformation group acting on $X$, and that $X$ is noncompact.
We show that
in this case $X\cong\RR^n$ is a real vector space, and that
$G$ splits as a semidirect product $G=G_x\ltimes\RR^n$. In view
of Lemma \ref{AbLem} and Proposition \ref{SimpleLem}, it
suffices to prove that $G^\circ$ cannot be simple if $X$ is
noncompact. We need some more facts about transformation groups.
For a group $H$ acting on a set $X$, we denote the fixed point set
by
\[
X^H=\{x\in X\mid H\cdot x=\{x\}\}.
\]
Recall the Malcev-Iwasawa Theorem, cp.~Salzmann \emph{et al.}
93.10.
\begin{Num}\textbf{\hspace{-1ex}Malcev-Iwasawa Theorem }
\label{MalcevIwasawa}
\em Let $G$ be a connected locally compact group. 

(1) There exists a maximal
compact subgroup $K\leq G$. If $L\leq G$ is a compact subgroup,
then $gLg^{-1}\leq K$ for a suitable $g\in G$.

(2) There exist closed subgroups $A_1,\ldots,A_m\leq G$
isomorphic to $(\RR,+)$, such that $G$ is directly spanned
as $G=A_1A_2\cdots A_mK$ (i.e. every $g\in G$ can in a unique way be
written as a product $g=a_1\cdots a_mk$, with $a_i\in A_i$ and $k\in K$).
In particular, $G$ is homeomorphic to
$K\times\RR^m$, and $G/K\cong\RR^m$ is contractible.
\end{Num}
We also need the following result.
\begin{Lem}
\label{pi3Lem}
Let $K$ be a compact connected Lie group. Then $K$ is
abelian if and only if $\pi_3(K)=0$.

\proof
It is a well-known result that $\pi_3(H)\cong\ZZ$ holds
for every almost simple compact Lie group $H$, see Onishchik
\cite{On} \S 17 Theorem 2, p~257.
A compact connected Lie group is topologically
a product of compact almost simple Lie groups and a torus; it is
abelian if and only if no nontrivial simple factors exist.
\qed
\end{Lem}

\begin{Cor}
\label{pi3Simple}
Let $G$ be an almost simple Lie group. If $G$ is not locally
isomorphic to $\SL_2\RR$, then $\pi_3(G)\neq 0$.

\proof
Since we consider only the third homotopy group, we may assume that
$G$ is centerless. Let $K\leq G$ be a maximal compact subgroup.
A direct inspection of Table V in Helgason \cite{He} Ch.~X,
p.~518 or Onishchik-Vinberg \cite{OnVin}
Table 9 p.~312 shows that $K$ is not abelian,
provided that $G\neq\mathrm{PSL}_2\RR$. Since we have a homotopy
equivalence $K\simeq G$ by the Malcev-Iwasawa
Theorem~\ref{MalcevIwasawa}, the claim follows.
\qed
\end{Cor}
For a proof of Whitehead's Theorem below see Spanier \cite{Spanier} 
Ch.~7 Sec.~6 Corollary 24.
\begin{Thm}
\label{Whitehead}
Let $f:X\too Y$ be a continuous map. If $X$ and $Y$ have the
homotopy types of CW-complexes (eg. if $X$ and $Y$ are manifolds,
see Bredon \cite{Bredon} V 1.6)
and if the induced map
$\pi(f)_\bullet:\pi_\bullet(X)\too\pi_\bullet(Y)$
is an isomorphism, then $f$ is a homotopy equivalence.
\end{Thm}
Combining these two results, we obtain the next lemma.

\begin{Lem}
\label{MaxHomEq}
Let $Z=H/H_z$ be a homogeneous space of a Lie group $H$. Assume that
$H$ and $H_z$ are connected. Let $K\leq H$ be a maximal compact subgroup,
such that $K_z$ is a maximal compact subgroup of $H_z$. Then
there is a homotopy equivalence $K/K_z\simeq H/H_z$.

\proof
The spaces $H/K$ and $H_z/K_z$ are contractible by the Malcev-Iwasawa
Theorem \ref{MalcevIwasawa}.
A fibration of CW-complexes with contractible fibres
(with a contractible base, resp.) induces a homotopy equivalence between
the base and the total space (the fibre and the total space, resp.)
by Whitehead's Theorem \ref{Whitehead}. Consider the fibrations
\[
\xymatrix@!{
& K/K_z \ar[d]^f\ar[dr]^{gf} \\
H_z/K_z \ar[r] & H/K_z \ar[r]^g\ar[d]  & H/H_z \\
& H/K}
\]
The maps $f$ and $g$ are homotopy equivalences, and so is
their composite $gf$.
\qed
\end{Lem}
The following theorem can be proved under much weaker assumptions;
however, the present version, which is taken from Bredon's book
\cite{Bredon} suffices for our purposes.

\begin{Thm}
\label{Cod1}
Let $X$ be a connected noncompact $n$-manifold, and let $K$ be a
compact Lie group acting smoothly or locally smoothly on $X$.
Assume that $K$  has an $(n-1)$-dimensional orbit. 
Then there are only the following two possibilities.

(a) $X/K\cong\RR$ and every $K$-orbit is principal. In this case,
$X$ decomposes as $X\cong K/K_x\times\RR$, with trivial $K$-action on $\RR$.

(b) $X/K\cong [0,\infty)$, and there exists precisely one nonprincipal
orbit $K\cdot z=Z\cong K/K_z$.
Let $K_x\leq K_z$ be the stabilizer of a 
point $x\in X\setminus Z$. Then $K_z/K_x\cong\SS^m$ is a sphere,
and $X$ is equivariantly homeomorphic to the $(m+1)$-vector bundle bundle
associated to the orthogonal sphere bundle $K_z/K_x\too K/K_x\too K/K_z$.

\proof
For the concept of a locally smooth action see Bredon's book
\cite{Bredon} Ch.~IV. A smooth action of a compact Lie group is
locally smooth by \emph{loc.cit.} VI 2.4; by IV 3.1,
principal orbits exist. Our claim is thus a slight reformulation
of \emph{loc.cit.} IV 8.1.
\qed
\end{Thm}

\begin{Cor}
\label{Cod1Cor}
If $n\geq 2$ and if $K$ has a fixed point, then $X$ is equivariantly
homeomorphic to euclidean space $\RR^n$, and $K$ acts transitively on
$\SS^{n-1}$.

\proof In this case, $Z=\{z\}$. A vector bundle over $Z$ is the same
as a real vector space $\RR^n$.
\qed
\end{Cor}
The proof of the next lemma was kindly pointed out by R. L\"owen.

\begin{Lem}
\label{DiscC}
Let $M$ be a connected manifold of dimension $n=\dim(M)\geq 3$. If
$D\subseteq M$ is a closed discrete subset, then the inclusion
$\xymatrix@1{M\setminus D\ \ar@{^(->}[r]&M}$
induces an isomorphism of fundamental groups.

\proof
Assume first that $D=\{x_0\}$ is a singleton. Let
$\phi:(\RR^n,0)\too(M,x_0)$ be a coordinate chart and put
$B_\eps=\{\phi(x)\mid |x|\leq\eps\}$. Then $M=B_2\cup(M\setminus B_1)$,
and $M\setminus B_1\simeq M\setminus\{x_0\}$. Moreover,
$B_2\cap (M\setminus B_1)=B_2\setminus B_1\simeq\SS^{n-1}$
is 1-connected, and
we can apply the Seifert-Van Kampen Theorem to the diagram
\[
{\xymatrix@!=2em{
& {B_2\setminus B_1} \ar[dl] \ar[dr] \\
B_2\ar[dr] && {M\setminus B_1}\ar[dl] \\
& M }}
{\xymatrix@!R=.4em{\\{} \ar@/^/[rrr]^{\text{apply $\pi_1(-)$}} &&&{}}}
{\xymatrix@!=2em{
& 1 \ar[dl]\ar[dr] \\ 
1 \ar[dr] && {\pi_1(M\setminus B_1)} \ar[dl] \\
& {\pi_1(M)} }}
\]
to conclude that $\pi_1(M\setminus\{x_0\})\too\pi_1(M)$ is an
isomorphism (see for example Bredon \cite{BredonTop} Ch.~III Cor.~9.5).
An easy induction shows that
$\pi_1(M\setminus D)\too\pi_1(M)$ is an isomorphism, provided that
$D$ is finite.

Assume now that $D\subseteq M$ is an arbitrary closed
discrete set not
containing the base point. We show first that
$\pi_1(M\setminus D)\too\pi_1(M)$ is surjective.
So let $\beta:[0,1]\too M$ be a path representing an element
$[\beta]_{\pi_1(M)}$
of $\pi_1(M)$, and let $U$ be a relatively compact open neighborhood
of the image of $\beta$. Then $U\cap D$ is finite, and by our
previous discussion, $\pi_1(U\setminus D)\too\pi_1(U)$ is
an isomorphism. Let $\alpha:[0,1]\too U\setminus D$ be a path
representing an element $[\alpha]_{\pi_1(U\setminus D)}$
in $\pi_1(U\setminus D)$ whose image $[\alpha]_{\pi_1(U)}$ in
$\pi_1(U)$ represents the same element $[\beta]_{\pi_1(U)}$ as $\beta$.
Then $\alpha$ and $\beta$ represent the same element in $\pi_1(M)$,
i.e. $[\alpha]_{\pi_1(M)}=[\beta]_{\pi_1(M)}$.

Next we prove that the map is injective. Let $\gamma$ be a path
representing an element of $\pi_1(M\setminus D)$, and assume that
its image $[\gamma]_{\pi_1(M)}$ is homotopic to the constant path.
Fix such a homotopy $F:[0,1]\times[0,1]\too M$, and let
$V$ be a relatively compact open neighborhood of the image of $F$.
Then $\gamma$ represents the identity element in $\pi_1(V)$.
By our previous discussion, $[\gamma]_{\pi_1(V\setminus D)}$
is also the identity element, and thus $[\gamma]_{\pi_1(M\setminus D)}$
is the identity element.
\qed
\end{Lem}

\begin{Def}
Let $G$ be a connected Lie transformation group acting transitively on a
manifold $X$. Assume that the following holds: for
every $x\in X$, the fixed point set
$X^{G_x}=\left\{y\in X\mid G_x\cdot y=\{y\}\right\}$ of
$G_x$ discrete (and necessarily closed),
and $G_x$ acts transitively on its complement $X\setminus X^{G_x}$.
Then we call $G$ \emph{almost 2-transitive} on $X$.
\end{Def}

\begin{Prop}
\label{1a2trs}
Let $G$ be a connected Lie group acting almost
2-transitively on a noncompact 1-connected $n$-manifold $X$, with
$n\geq 3$. Then $G$ acts 2-transitively on $X$, and $X\cong\RR^n$ is
contractible.

\proof
By assumption, $\pi_0(G)=1=\pi_1(X)$. From the exact sequence 
\[
\xymatrix{
\pi_1(X)\ar[r] & \pi_0(G_x)\ar[r] & \pi_0(G)}
\]
we see that $G_x$ is connected, and $X\setminus X^{G_x}$ is
1-connected by Lemma \ref{DiscC}. Let $y\in X\setminus X^{G_x}$,
so $G_x/G_{x,y}\cong X\setminus X^{G_x}$. The exact sequence
\[
\xymatrix{
\pi_1(X\setminus X^{G_x})\ar[r]&\pi_0(G_{x,y})\ar[r]&\pi_0(G_x)}
\]
shows that $G_{x,y}$ is connected.
We choose a maximal compact subgroup $K\leq G$ such that
$K_x\leq G_x$  and $K_{x,y}\leq G_{x,y}$ are maximal compact
subgroups. Note that each of these groups is connected.

Let $\FF_2$ denote the field with
two elements. Then $H_k(X;\FF_2)=0=H_k(X\setminus X^{G_x};\FF_2)$ for all
$k\geq n$, because $X$ and $X\setminus X^{G_x}$ are noncompact manifolds,
see Bredon \cite{BredonTop} Ch.~VI Cor.~7.12 e.g.
Moreover, $\dim(H_n(X,X\setminus X^{G_x};\FF_2))=|X^{G_x}|$
(by excision), and
the exact sequence
\[
\xymatrix{{}\ar[r]&
H_n(X;\FF_2)\ar[r] &
H_n(X,X\setminus X^{G_x};\FF_2)\ar[r]&
H_{n-1}(X\setminus X^{G_x};\FF_2)\ar[r]&{}}
\]
shows that $\dim(H_{n-1}(X\setminus X^{G_x};\FF_2))\geq |X^{G_x}|$.
But $X\setminus X^{G_x}\simeq K_x/K_{x,y}$ is homotopy equivalent to
a compact connected manifold of dimension strictly less than $n$
by Lemma \ref{MaxHomEq}. Thus
$\dim(K_x/K_{x,y})=n-1$, and $H_{n-1}(K_x/K_{x,y};\FF_2)\cong\FF_2$.
In particular, $|X^{G_x}|=1$, and thus $G_x$ acts transitively on
$X\setminus\{x\}$, i.e. $G$ is 2-transitive.
The orbit $K_x\cdot y\cong K_x/K_{x,y}$
has codimension 1 in $X$. By Corollary \ref{Cod1Cor},
$X\cong\RR^n$.
\qed
\end{Prop}

\begin{Prop}
\label{Contractible}
Let $G$ be a 2-transitive Lie group acting on a noncompact
connected manifold $X$ of dimension $\dim(X)=n\geq 3$.
Then $X$ is 1-connected, and thus by Proposition \ref{1a2trs}
$X\cong\RR^n$ is contractible.

\proof
Replacing $G$ by the universal cover of its connected component,
we may assume that $G$ is
1-connected. Thus we have an isomorphism $\pi_1(X)\cong\pi_0(G_x)$.
Let $H=(G_x)^\circ$. The exact sequence
\[
\xymatrix{{}\ar[r] &
\pi_1(G)\ar[r]&\pi_1(G/H)\ar[r]&\pi_0(H)\ar[r]&{}}
\]
shows that $Z=G/H$ is 1-connected, and
$\xymatrix{Z=G/H\ar[r]^p&G/G_x=X}$ is the universal covering of $X$.
Let $F=p^{-1}(x)$. We claim that $H$ acts transitively on $Z\setminus F$.
Let $z\in Z\setminus F$. Then $p(H\cdot z)=H\cdot p(z)=X\setminus \{x\}$,
because $H$ is open in $G_x$ and because $X\setminus\{x\}$ is connected.
Thus the $H$-orbits in $Z\setminus F$ are $n$-dimensional and hence
open. But $Z\setminus F$ is connected, so $H\cdot z=Z\setminus F$.
Thus $G$ acts almost 2-transitively on $G/H$. By Proposition \ref{1a2trs},
the action is in fact 2-transitive and thus $H=G_x$.
\qed
\end{Prop}
Finally, we consider the low-dimensional cases.

\begin{Lem}
\label{LowDim}
Let $G$ be a 2-transitive Lie group acting on a noncompact
connected manifold $X$ of dimension $\dim(X)=n\leq 2$.
Then $X\cong\RR^n$ is contractible.

\proof
The only noncompact connected 1-manifold is $\RR$.

If $\dim(X)=2$, then $H=G^\circ$ acts 2-transitively by Lemma \ref{CC2trs}.
Assume that $X$ is a noncompact surface. If $H$ is not simple, then
$X\cong\RR^2$ by Lemma \ref{AbLem}. Assume that $H$ is simple.
Since $H$ acts 2-transitively, $\dim(H)\geq 4$, and in particular
$H\neq\mathrm{PSL}_2\RR$. Let $K$ be a maximal compact subgroup of $H$.
Since $H\neq\mathrm{PSL}_2\RR$, the group $K$ is not abelian, and
thus $\pi_3(H)=\pi_3(K)\neq 0$ by Corollary \ref{pi3Simple}.
Since $X$ is a noncompact surface, the only homotopy group of
$X$ which is possibly nontrivial is the fundamental group.
The exact sequence
\[
\xymatrix{{}\ar[r]&
\pi_4(X)\ar[r]&
\pi_3(H_x)\ar[r]&
\pi_3(H)\ar[r]&
\pi_3(X)\ar[r]&{}}
\]
shows that $\pi_3(H_x)\neq 0$, hence $H_x$ contains a connected
nonabelian compact subgroup $L\leq (H_x)^\circ$ by Lemma
\ref{pi3Lem}. Choose $y\in X\setminus X^L$.
Then $L\cdot y\cong\SS^1$ is a circle,
and $X\cong\RR^2$ by Corollary \ref{Cod1Cor}.
(We will see below that this case is impossible.)
\qed
\end{Lem}

\begin{Cor}
If $G$ is an effective 2-transitive Lie group acting on $X$ and
if $X$ is noncompact, then $G^\circ$ is not simple.

\proof
By Proposition \ref{Contractible} and Lemma \ref{LowDim}, the space
$X$ is contractible. Assume that $H=G^\circ$ is simple. Since
$X$ is contractible, the stabilizer $H_x$ is connected, and
$\xymatrix@1{H_x\ \ar@{^(->}[r]&H}$ is a homotopy equivalence by
Whitehead's Theorem \ref{Whitehead}. Thus $H_x$
contains a maximal compact subgroup $K$ of $G$. But the maximal
compact subgroups in a noncompact simple Lie group are maximal
subgroups, see Helgason \cite{He} Ch.~VI Ex.~A3(iv)
p.~276 and p.~567,
hence $H_x=K$, so $X=H/K$ is a
Riemannian symmetric space, and $H$ preserves the metric;
in particular, $H$ has infinitely many orbits on $X\times X$,
contradicting Theorem \ref{2-TrsLie}.
\qed
\end{Cor}
Combining this result with our previous analysis, we obtain the
following final result for the noncompact case.

\begin{Thm}
Let $G$ be a locally compact and $\sigma$-compact group acting
effectively and 2-transitively on a locally compact, noncompact, not
totally disconnected space $X$. Then $G$ is a Lie group and 
has a normal vector subgroup $\RR^n\unlhd G$ acting regularly
on $X\cong\RR^n$. Moreover,
$G=G_x\ltimes\RR^n$ is a semidirect product of
$\RR^n$ and a point stabilizer $G_x$.

If $n=1$, then $G=\mathrm{AGL}_1\RR=\GL_1\RR\ltimes\RR$.

If $n\geq 2$, then the connected stabilizer $(G_x)^\circ$
is a reductive linear Lie group acting transitively on the nonzero
vectors of $\RR^n$.
All possibilities for $(G_x)^\circ$ (and thus for
$G^\circ=(G_x)^\circ\ltimes\RR^n$)
are determined in Theorem \ref{TrsLinGrp}.
None of these actions is 3-transitive.
If $G$ acts sharply 2-transitive
(i.e. if the two-point stabilizers are trivial),
then $n=1,2,4$, and $G_x=(G_x)^\circ$
is one of the groups given in Corollary \ref{Sh2Trs}.

If $X$ is a complex manifold and if $G$ preserves the complex
structure, then $G^\circ$ is one of the groups given in
Theorem \ref{TrsLinGrp} (d).
\qed
\end{Thm}
In particular, we obtain Knop's result \cite{Knop}
for characteristic 0:

\begin{Cor}
If $G$ is a complex connected Lie group which acts complex
analytically
and 2-tran\-si\-ti\-ve\-ly on a complex manifold $X$, then either 

(a) $X\cong\CC\mathrm{P}^n$ and $G=\mathrm{PSL}_{n+1}\CC$, or

(b) $X\cong\CC^n$ and $G_x=\SL_n\CC$ or $G_x=\GL_n\CC$, or

(c) $X\cong\CC^{2n}$ and $G_x\cong\Sp_{2n}\CC$ or
$G_x\cong\Sp_{2n}\CC\cdot\CC^*$.

\noindent
In case (b) or (c), $G\cong G_x\ltimes\CC^n$.
\end{Cor}
Note however that in the case $X\cong\CC^n$, there exist several
real, noncomplex 2-transitive groups by Theorem \ref{TrsLinGrp}.

\section{Transitive groups acting on $\RR^m\setminus 0$}

In this section we classify all closed linear Lie groups
acting transitively on the nonzero vectors of a finite
dimensional real vector space. The main ingredients in
our proof are the classification of compact connected linear
Lie groups acting transitively on spheres, and the representation
theory of compact Lie groups.

\begin{Thm}
\label{TrsSphere}
Let $m\geq 2$, and let
$K\leq\SO(m)$ be a closed connected subgroup. If $K$ acts
transitively on the sphere $\SS^{m-1}=\{x\in\RR^m\mid |x|=1\}$,
then $K$ is (up to automorphisms of $\SO(m)$)
one of the following groups.
\begin{center}
\begin{tabular}{l|ll}
$K$ & $\RR^m$ & $m$ \\ \hline
$\SO(n)$ & $\RR^n$ & $n$ \\
$\SU(n)$ & $\CC^n$ & $2n$ \\
$\U(n)$  & $\CC^n$ & $2n$ \\
$\Sp(n)$ & $\HH^n$ & $4n$ \\
$\Sp(n)\cdot\U(1)$  & $\HH^n$ & $4n$ \\
$\Sp(n)\cdot\Sp(1)$ & $\HH^n$ & $4n$ \\
$\mathrm{G}_2$ & $\mathrm{Pu}(\mathbb{O})$ & $7$ \\
$\Spin(7)$ & $\mathbb{O}$ & $8$ \\
$\Spin(9)$ & $\mathbb{O}\oplus\mathbb{O}$ & $16$
\end{tabular}
\end{center}
Besides the standard inclusions between the classical groups
(eg. $\Sp(k)\leq\SU(2k)\leq\U(2k)\leq\SO(4k)$ etc.), the
inclusions in the special dimensions $7$, $8$, and $16$ are as follows
(note that $\SU(4)\cong\Spin(6)$ and $\Sp(2)\cong\Spin(5)$).
\[
\xymatrix@R=2em@C=2em{
7 & 8        &&& 16 \\
\SO(7) \ar@{-}[d]
& \SO(8) \ar@{-}[d]\ar@{-}[dr]  &&& \SO(16) \ar@{-}[d]\ar@{-}[dl] \\
\mathrm{G}_2 & \Spin(7) \ar@{-}[d] & \U(4) \ar@{-}[d]\ar@{-}[dl] & \Spin(9)
& \U(8) \ar@{-}[d]\ar@{-}[dr] \\
& \SU(4)   \ar@{-}[d] & \U(1)\cdot\Sp(2) \ar@{-}[dl] &&
\U(1)\cdot\Sp(4) \ar@{-}[dr] & \SU(8) \ar@{-}[d] \\
& \Sp(2)   &&&& \Sp(4)
}
\]
All exceptional phenomena in dimensions $7$, $8$, $16$ are related
to the Cayley plane $\mathrm{PG}_2\OO$;
see Salzmann \emph{et al.} \cite{CPP} Ch.~I.

\proof
We could use the classification of compact Lie group acting
transitively on spheres, due to Montgomery-Samelson \cite{MoSa}
and Borel \cite{BoS1,BoS2}; see also Poncet \cite{Poncet},
Besse \cite{Besse} 7.13,
Onishchik \cite{On} \S 18 Theorem 3 (i),
Salzmann \emph{et al.} \cite{CPP} 96.20--23,
and Kramer \cite{KrMem} Ch.~6.
However, we actually need only the classification of transitive
subgroups of $\SO(n)$, which is much easier and follows directly
from Onishchik's classification of factorizations of $\SO(n)$,
see \cite{On} p.~227 and \cite{GOV} II \S4.5, p.~144.
\qed
\end{Thm}

\begin{Cor}
\label{ComTrsS}
Suppose that $m\geq 3$. Then the commutator group $K'=[K,K]$
acts also transitively on $\SS^{m-1}$.

\proof
This follows either by direct inspection of the list, or by a simple
homotopy-theoretic argument, see Onishchik \cite{On} \S 5 Prop.~9
or Grundh\"ofer-Knarr-Kramer \cite{GKK} Lem.~1.3.
\qed
\end{Cor}

\begin{Prop}
\label{TrsProp1}
Let $V$ be a finite dimensional real vector space of dimension
$m\geq 3$, and let
$H\leq\GL(V)$ be a closed connected group. Suppose that
$H$ acts transitively on the nonzero vectors of $V$. Then
$H$ is reductive, and consequently $H'=[H,H]$ is a semisimple
subgroup of $\SL(V)$. If $K\leq H$ is a maximal compact subgroup
and if $|-|$ is a $K$-invariant norm on $V$, then
$K$ and $[K,K]$ act transitively on the sphere
$\SS^{m-1}=\{x\in V\mid |x|=1\}$.
Consequently, $K$ is one of the groups given in Theorem \ref{TrsSphere}.

\proof
Clearly, $H$ acts irreducibly on $\RR^m$, and thus
$H$ is reductive, see Salzmann \emph{et al.} \cite{CPP} 95.2, 95.6.
Let $K\leq H$ be a maximal compact subgroup, let
$|-|$ be a $K$-invariant norm, and let $v$ be a vector with $|v|=1$.
There is a homotopy equivalence
$K/K_v\simeq H/H_v\cong\RR^m\setminus\{0\}$ by Lemma \ref{MaxHomEq}.
Thus $K/K_v\simeq\SS^{m-1}$, and in particular $\dim(K/K_v)=m-1$.
By domain invariance, $K\cdot x=\SS^{m-1}$, and by
Corollary \ref{ComTrsS}, the commutator group $[K,K]$ is also transitive
on $\SS^{m-1}$.
\qed
\end{Prop}

\begin{Cor}
There exists a simple subgroup $H_1\leq H$ such that $K\cap H_1$
acts transitively on $\SS^{m-1}$.

\proof
By Proposition \ref{TrsProp1},
$[K,K]\leq[H,H]$ acts transitively on $\SS^m$. Theorem
\ref{TrsSphere} shows in each case that $K$ has a simple factor
which acts transitively on $\SS^{m-1}$. We can choose a
simple factor $H_1\leq [H,H]$ containing this simple factor of
$K$.
\qed
\end{Cor}
Now we determine all possibilities for the semisimple group $[H,H]$.
The following facts will be useful.
\begin{Thm}
\label{DoubleCartan}
Let $H\leq\SL_m\RR$ be a closed semisimple subgroup. After conjugation
with some element in $\SL_m\RR$, the Cartan decomposition of 
$\fsl_m\RR=\fso(m)\oplus\mathfrak{p}$ is also a Cartan decomposition of
$\fh$,
\[
\fh=(\fh\cap\fso(m))\oplus(\fh\cap\mathfrak{p}).
\]

\proof
See Onishchik-Vinberg \cite{OnVin} Ch.~5 Theorem 4, p.~261.
\qed
\end{Thm}

\begin{Lem}
\label{PGenerates}
Let $\fg=\fk\oplus\fp$ be a Cartan decomposition of a simple noncompact
Lie algebra. Then 
\[
[\fp,\fp]=\fk.
\]
In particular, if $\mathfrak{l}\leq\fg$ is a subalgebra containing
$\fp$, then $\mathfrak{l}=\fg$.

\proof Let $\fh=[\fp,\fp]\leq\fk$. Using the Jacobi identity, it is
clear that $[\fk,\fh]\leq\fh$. Since
$[\fp,\fh]\leq\fp$, we conclude that $\fh+\fp$ is an ideal in $\fg$,
hence $\fh+\fp=\fg=\fk+\fp$. This is a direct sum decomposition
with $\fh\leq\fk$, therefore $\fh=\fk$.
\qed
\end{Lem}
We use the following facts from representation theory. We rely on
Tits \cite{Ti2} and \cite{Ti3}, on
the Reference Chapter in Onishchik-Vinberg \cite{OnVin}, and on
Chapter 9 in Salzmann \emph{et al.} \cite{CPP}.

\begin{Num}\textbf{Remarks on representation theory}
Associated to a real semisimple Lie group $G$ is its weight lattice;
every weight $\pi$ is a linear combination with integral coefficients
of the so-called fundamental
weights $\pi_i$. The dominant weights (i.e. the
weights where all coefficients are nonnegative) correspond
bijectively to the equivalence classes of (finite dimensional)
complex irreducible
representations of $G$. We denote the complex irreducible module
associated to the dominant weight $\pi$ by $R(\pi)$. The Galois group of
$\CC/\RR$ acts on the weight lattice via $\pi\mapstoo\bar\pi$;
if $\pi=\bar\pi$ then
$R(\pi)$ admits a real or a quaternionic structure (i.e. a complex
semilinear endomorphism $\phi$ with $\phi^2=\id$ resp.~$\phi^2=-\id$).
There is a corresponding map $\beta$ from the invariant weights into the
Brauer group $\mathrm{Br}(\RR)=\{\RR,\HH\}$. There are three types
of real irreducible representations of $G$:

(1) If $\pi\neq\bar\pi$, then $\pi$ is of complex type and
$_\RR R(\pi)=R(\pi)$ is a real irreducible $G$-module.

(2) If $\pi=\bar\pi$ and $\beta_\pi=\RR$, then $\pi$ is
of real type. There exists a unique real irreducible $G$-module
$_\RR R(\pi)$, and $R(\pi)\cong{}_\RR R(\pi)\otimes\CC$.

(3) If $\pi=\bar\pi$ and $\beta_\pi=\HH$, then $\pi$
is of quaternionic type. In this case, $_\RR R(\pi)=R(\pi)$
is a real irreducible module, and $R(\pi)$ admits a quaternionic
structure.
\end{Num}
The Reference Chapter in Onishchik-Vinberg \cite{OnVin} contains tables
which we will frequently use to explain how exterior powers
$\bigwedge\nolimits^k V$, symmetric powers $S^kV$, and
tensor products $V^{\otimes k}$ of complex irreducible modules $V$
can be decomposed.

We fix some notation. We denote the ring of $n\times n$-matrices with
entries in $F$ by $F(n)$. The notation for Lie groups and algebras
is as in Onishchik-Vinberg \cite{OnVin}. In particular, we denote the
quaternion unitary group $\U_n\HH$ by $\Sp(n)$. Let
\[
S(n)=\{X\in\RR(n)\mid X^T=X\}\quad\text{ and }\quad
S_0(n)=\{X\in S(n)\mid \mathrm{tr}(X)=0\}
\]
The Cartan decomposition of $\fsl_m\RR$ is $\fsl_m\RR=\fso(m)\oplus S_0(m)$.
\begin{Lem}
\label{SymLemma}
Let $H\leq\SO(n)$ be a subgroup. Then there is a standard
action of $H$ on $S_0(n)$ which is given by
$X\mapstoo hXh^T$. The following actions of subgroups
of $\SO(n)$ on $S_0(n)$ are $\RR$-irreducible for $n\geq 2$:
\begin{center}
\begin{tabular}{l|l}
$H$ & $n$ \\ \hline
$\SO(n)$ & $n$ \\
$\mathrm{G}_2$ & $7$ \\
$\Spin(7)$ & $8$
\end{tabular}
\end{center}

\proof
For the orthogonal groups, this follows from the fact that every
quadratic form can be diagonalized by conjugation with an element of
$\SO(n)$. If $0\neq U\leq S_0(n)$ is an $\SO(n)$-invariant subspace,
then $U$ intersects the space $A$ consisting of diagonal traceless
matrices nontrivially. Moreover, the $\SO(n)$-stabilizer of $A$
induces the symmetric group $\mathrm{Sym}(n)$ on the diagonal
matrices. Thus $U$ contains $A$, whence $U=S_0(n)$.
(The action of $\SO(n)$ is in fact a \emph{polar action}; it is
the isotropy representation of the Riemannian symmetric space
$\SL_n\RR/\SO(n)$. The principal orbits are isoparametric
submanifolds; the orbit types form the building $\Delta(\SL_n\RR)$,
i.e. the real projective geometry of rank $n-1$).

For the groups $\mathrm{G}_2$ and $\Spin(7)$ we use representation
theory.
The complex module corresponding to the action of $\mathrm{G}_2$ on
the set $\mathrm{Pu}(\OO)\cong\RR^7$ of pure octonions is $R(\pi_1)$, and 
\[
S(7)\otimes\CC\cong S^2R(\pi_1)=R(2\pi_1)\oplus\CC.
\]
The Galois group of $\CC/\RR$ act trivially on the weights, and 
all representations of $\mathrm{G}_2$ are of real type, see
Tits \cite{Ti2} p.~42. Thus $S_0(7)\cong {}_\RR R(2\pi_1)$ is a
real irreducible $\mathrm{G}_2$-module.

The reasoning for $\Spin(7)$ is similar. The complex module for the
$8$-dimensional representation is $R(\pi_3)$, and
$S^2R(\pi_3)=R(2\pi_3)\oplus\CC$, see Onishchik-Vinberg \cite{OnVin} p.~301.
This module is of real type, see Tits \cite{Ti2} p.~31, hence
${}_\RR R(2\pi_3)\otimes\CC\cong R(2\pi_3)$.
Thus $S_0(8)\cong{}_\RR R(2\pi_3)$ is a real irreducible
$\Spin(7)$-module.
\qed
\end{Lem}
Now we determine all closed noncompact semisimple subgroups of $\SL_m\RR$
which contain one of the compact groups in Theorem \ref{TrsSphere}.
\begin{Prop}
\label{ContainsSO(n)}
Let $H\leq \SL_m\RR$ be a closed noncompact semisimple subgroup, for
$m\geq 3$. If $\SO(m)<H$, then $H=\SL_m\RR$. If $m=7$ and if
$\mathrm{G}_2<H$, then $H=\SL_7\RR$. If $m=8$ and if $\Spin(7)<H$,
then $H=\SL_8\RR$.

\proof
Let $\fh$ denote the Lie algebra of $H$. We may assume that
$\fh\leq\fsl_m\RR=\fso(n)\oplus S_0(n)$ is embedded as in
Theorem \ref{DoubleCartan}. Let $K=\SO(m)\cap H$.
Then $\fh\cap S_0(m)$ is a nonzero $K$-module. Since we assume
that $K=\SO(m)$ (resp. that $\mathrm{G}_2\leq K$ or that $\Spin(7)\leq K$),
it follows from Lemma \ref{SymLemma} that $\fh\cap S_0(m)=S_0(m)$, and
thus $\fh=\fsl_m\RR$ by Lemma \ref{PGenerates}.
\qed
\end{Prop}
The case of $\Spin(9)$ acting on $\OO\oplus\OO$ is
different. Consider the affine Cayley plane $\mathrm{AG}_2\OO$,
and let $H\leq\GL_{16}\RR$ denote the stabilizer of the origin
in its collineation group.
Then $H\cong(\Spin_{9,1}\RR)^\circ$, and $H\cap\mathrm{F}_4\cong\Spin(9)$,
see Salzmann \emph{et al.} \cite{CPP}
Sec.~15, in particular 15.6, and p.~628.

\begin{Prop}
\label{ContainsSpin(9)}
Let $H\leq\SL_{16}\RR$ be a closed noncompact semisimple subgroup with
$\Spin(9)<H$. Then $H$ is one of the groups
$(\Spin_{9,1}\RR)^\circ$ or $\SL_{16}\RR$.

\proof
The idea is to decompose the $\Spin(9)$-module
$S_0(16)$ into irreducible submodules.
We denote the Lie algebra of $\Spin(9)\leq\SO(16)$
by $\fspin(9)\leq\fso(16)$.
Let $V=\OO\oplus\OO$ denote the 'natural' $\Spin(9)$-module.
Then $S(16)\cong S^2V$.

The Galois group of $\CC/\RR$ acts trivially on all weights of
$\Spin(9)$, and all
fundamental representations of $\Spin(9)$ are of real type,
see Tits \cite{Ti2} p.~31. We label the fundamental weights of
$\Spin(9)$ in the standard way (as in Onishchik-Vinberg \cite{OnVin}
p.~293 or in Tits \cite{Ti2} p.~30). The complex module
$R(\pi_4)$ associated to the weight $\pi_4$ is the $\Spin$-module,
$R(\pi_4)\cong V\otimes\CC$. By Onishchik-Vinberg \cite{OnVin} p.~301,
we have a decomposition
\[
S^2R(\pi_4)\cong R(2\pi_4)\oplus R(\pi_1)\oplus\CC 
\]
Therefore we obtain a decomposition into real irreducible modules
\[
S^2V\cong{}_\RR R(2\pi_4)\oplus{}_\RR R(\pi_1)\oplus\RR
\quad\text { and }\quad
S_0(16)\cong{}_\RR R(2\pi_4)\oplus{}_\RR R(\pi_1)
\]
The dimensions of the modules $_\RR R(\pi_1)\cong\RR^9$,
and $_\RR R(2\pi_4)$ are $9$ and $126$, respectively, see
Onishchik-Vinberg \cite{OnVin} p.~301.

Now suppose that $H\leq\SL_{16}\RR$ is embedded as in
Theorem \ref{DoubleCartan}. Since $\Spin(9)\leq H$, the Lie
algebra $\fh$ of $H$ is a $\Spin(9)$-module. Consider the nonzero
$\Spin(9)$-module $M=\fh\cap S_0(16)$. If $M=S_0(16)$, then
$\fh=\fsl_{16}\RR$ by Lemma \ref{PGenerates}.
Since we know that $\fso_{9,1}\RR\leq\fsl_{16}\RR$, the case
$M={}_\RR R(\pi_1)$ is possible, and $[M,M]=\fspin(9)$ in this
case, so $\fh=\fso_{9,1}\RR$ is uniquely determined. Finally,
suppose that $M={}_\RR R(2\pi_4)$. Then $\dim(\fspin(9)+M)=162$;
no such semisimple Lie algebra  (with $\fspin(9)$ as maximal
compact subalgebra) exists.
\qed
\end{Prop}
Now we consider semisimple groups which contain $\SU(n)$.
Let $H(n)=\{X\in\CC(n)\mid \bar X^T=X\}$ denote the set of
hermitian matrices, and let $H_0(n)=\{X\in H(n)\mid\mathrm{tr}(X)=0\}$.
Note that $\fsu(n)=\ti H_0(n)$ and $\fu(n)=\ti H(n)$. 

Let $\Theta=\scriptsize\begin{pmatrix}1 & \\ & -1\end{pmatrix}$
and $\I=\scriptsize\begin{pmatrix} & 1 \\ -1 & \end{pmatrix}$.
We identify $\CC$ with the subalgebra of $\RR(2)$ spanned
by $1,\I$; this subalgebra has a vector space complement
$\CC\Theta$ which is spanned by $\Theta$ and $\I\Theta$.
Every matrix $X\in\RR(n)\otimes\RR(2)\cong\RR(2n)$ decomposes
uniquely as
\[
X=X_1\otimes 1+X_2\otimes\I+(X_3\otimes 1+X_4\otimes\I)\Theta.
\]
We call $X_1\otimes1+ X_2\otimes\I$ the \emph{complex part} of $X$ and
$X_3\otimes\Theta+X_4\otimes\I\Theta$ the \emph{anti-complex part} of $X$.
There is a natural injection
$\xymatrix@1{\CC(n)\ \ar@{^(->}[r]&\RR(2n)}$ which is
given
by $A+B\ti\mapstoo A\otimes 1+B\otimes\I$, for $A,B\in\RR(n)$,
and $A+B\ti\in H(n)$ if and only if $A\in S(n)$ and
$B\in\fso(n)$. More generally, we have
\[
(X_1\otimes 1+ X_2\otimes\I+X_3\otimes\Theta+X_4\otimes\I\Theta)^T=
X_1^T\otimes1-X_2^T\otimes\I+X_3^T\otimes\Theta+X_4^T\otimes\I\Theta.
\]
Using this identity, it is not hard to show that
$S_0(2n)=H_0(n)\oplus S(n)\otimes\CC\Theta$.
The corresponding decompositions into $\SU(n)$-modules is
\[
S_0(2n)\cong\fsu(n)\oplus S^2V.
\]
where $V\cong\CC^n$ is the natural module (the symmetric 
power of $V$ being taken over $\CC$!).
Clearly, $\fsu(n)$ is an irreducible real $\SU(n)$-module for all
$n\geq 2$. By Onishchik-Vinberg \cite{OnVin} p.~300, the complex
module $S^2R(\pi_1)\cong R(2\pi_1)$ is irreducible for all $n\geq 2$;
and $2\pi_1$ is of complex type, $_\RR R(2\pi_1)=R(2\pi_1)$,
provided that $n\geq 3$. Thus $S^2V$ is a real irreducible $\SU(n)$-module
for $n\geq 3$.

\begin{Prop}
\label{ContainsSU(n)}
Let $H\leq \SL_{2n}\RR$ be a closed semisimple noncompact subgroup with
$\SU(n)<H$, for $n\geq 3$. Then $H$ is one of the groups
$\SL_n\CC$, $\Sp_{2n}\RR$, or $\SL_{2n}\RR$.

\proof
We assume that $\fh\leq \fsl_{2n}\RR$ is embedded as in Theorem
\ref{DoubleCartan}. Thus $\fh\cap S_0(2n)$ is a nonzero
$\SU(n)$-module. We showed above that
$S_0(2n)=H_0(n)\oplus S(n)\otimes\CC\Theta$
is a decomposition into real irreducible $\SU(n)$-modules.
The real dimensions of these two modules are $n^2-1$ and $n^2+n$,
respectively, so they are not isomorphic.
Moreover, the following subalgebras exist:
\[
\fsl_n\CC=\fsu(n)\oplus H_0(n) 
\quad\text{ and }\quad
\fsp_{2n}\RR=\fu(n)\oplus S(n)\otimes\CC\Theta.
\]
Since
\begin{gather*}
[H_0(n),H_0(n)]=\fsu(n),\qquad
[S(n)\otimes\CC\Theta,S(n)\otimes\CC\Theta]=\fu(n)\qquad\text{ and }\\
[S_0(2n),S_0(2n)]=\fso(2n)
\end{gather*}
by Lemma \ref{PGenerates}, the algebras
$\fsl_n\CC$, 
$\fsp_{2n}\RR$, and
$\fsl_{2n}\RR$
are the only possibilities for $\fh$.
\qed
\end{Prop}
Now we consider the case where $\fsp(n)\leq\fh$, for $n\geq 2$.
Let $V=\CC^{2n}=\HH^n$ denote the natural $\Sp(n)$-module.
The complex $\SU(2n)$-module corresponding to $\fsu(2n)$
is
\[
\fsu(n)\otimes\CC\cong\fsl_{2n}\CC=\{X\in\CC(2n)\mid\mathrm{tr}(X)=0\}.
\]
Thus we have to consider the complex $\Sp(n)$-module $V\otimes V$.
By Onishchik-Vinberg \cite{OnVin} p.~302,
$V\otimes V=R(\pi_1)\otimes R(\pi_1)\cong R(\pi_2)\oplus R(2\pi_1)
\oplus\CC$. The Galois group of $\CC/\RR$ acts trivially on the
fundamental weights, and $\pi_k$ is of real type if and only if
$k$ is even, see Tits \cite{Ti2} p.~34.
Thus $\pi_2$ and $2\pi_1$ are of real type.
Moreover, $R(2\pi_1)=\fsp_{2n}\CC$, and we obtain a
decomposition into real irreducible $\Sp(n)$-modules
\[
\fsu(2n)\cong\fsp(n)\oplus{}_\RR R(\pi_2).
\]
The dimension of $_\RR R(\pi_2)$ is $(2n+1)(n-1)$, and
$\dim(\fsp(n))=(2n+1)n$; in particular, the modules are not isomorphic.
The complex $\fsp(n)$-module $S(2n)\otimes\CC\cong S^2 R(\pi_1)$
decomposes as
\[
S^2R(\pi_1)\cong R(2\pi_1),
\]
see Onishchik-Vinberg \cite{OnVin} p.~302.
Thus the decomposition into real irreducible modules is
\[
S^2V\cong{}\fsp(n)\oplus\fsp(n)
\]
and
\[
S_0(4n)\cong\fsp(n)\oplus\fsp(n)\oplus\fsp(n)\oplus{}_\RR R(\pi_2).
\]
\begin{Prop}
\label{ContainsSp(n)}
Let $H\leq \SL_{4n}\RR$ be a closed noncompact semisimple subgroup with
$\Sp(n)<H$, for $n\geq 2$. Then $H$ is one of the groups
$\SL_n\HH$, $\SL_n\HH\cdot\Sp(1)$, $\SL_{2n}\CC$, $\Sp_{2n}\CC$,
$\SL_{4n}\RR$, or $\Sp_{4n}\RR$.

\proof 
We assume that $\fh\leq \fsl_{4n}\RR$ is embedded as in Theorem
\ref{DoubleCartan}. Thus $\fh\cap S_0(4n)$ is a nonzero
$\Sp(n)$-module. Note that the subalgebras
$\fsl_n\HH$, $\fsp_{2n}\CC$, $\fsl_{2n}\CC$ and $\fsp_{4n}\RR$ exist
in $\fsl_{4n}\RR$. We have
\begin{align*}
\fsl_n\HH\cap S_0(4n)      & \cong {}_\RR R(\pi_2) \\
\fsp_{2n}\CC\cap S_0(4n)   & \cong \fsp(n) \\
\fsl_{2n}\CC\cap S_0(4n)   & \cong \fsp(n)\oplus{}_\RR R(\pi_2) \\
\fsp_{4n}\RR\cap S_0(4n)   & \cong \fsp(n)\oplus\fsp(n) \\
\fsl_{4n}\RR\cap S_0(4N)   & \cong S_0(4n)
\end{align*}
Let $a\in\Sp(1)=\mathrm{Cen}_{\SO(4n)}\Sp(n)$ be a pure element,
$a+\bar a=0$. Then $a$ defines a real symplectic structure $\omega_a$
on $\RR^{4n}=\HH^n$ by $\omega_a(u,v)=\mathrm{Re}(\bar u^Tav)$.
Thus we see that all copies of $\fsp(n)\leq S_0(4n)$ are conjugate
under the group $\Sp(1)$.
In view of Lemma \ref{PGenerates}, we immediately
obtain the following results.

If $\fh\cap S_0(4n)={}_\RR R(2\pi_1)$, then
$\fsl_n\HH\leq\fh$ and thus $\fh=\fsl_n\HH$ or
$\fh=\fsl_n\HH\oplus\fsp(1)$.

If $\fh\cap S_0(4n)\cong\fsp(n)$, then $\fh$ in conjugate to
the algebra $\fsp_{2n}\CC$.

If $\fh\cap S_0(4n)\cong\fsp(n)\oplus{}_\RR R(2\pi_1)$, then
$\fh$ is conjugate to $\fsl_{2n}\CC$.

If $\fh\cap S_0(4n)\cong\fsp(n)\oplus\fsp(n)$, then
$\fh$ is conjugate to $\fh=\fsp_{4n}\RR$.

If $\fh\cap S_0(4n)= S_0(4n)$, then $\fh=\fsl_{4n}\RR$.

Thus, if $\fh\cap S_0(4n)$ contains
$\fsp(n)\oplus{}_\RR R(\pi_2)$ or $\fsp(n)\oplus\fsp(n)$, then
$\fsu(2n)\leq\fh$, and we showed
in Proposition \ref{ContainsSU(n)} that in this case either
$\fh=\fsp_{4n}\RR$ or $\fh=\fsl_{2n}\CC$.

There are no other cases to be considered.
\qed
\end{Prop}
Finally, we discuss the low dimensional cases.
\begin{Prop}
\label{ContainsSU(2)}
Let $H\leq\SL_4\RR$ be a closed semisimple noncompact Lie group
containing $\SU(2)$. Then $H$ is one of the groups
$\SL_2\CC$, $\Sp_4\RR$, or $\SL_4\RR$.

\proof
As we noted before Proposition \ref{ContainsSU(n)}, we have
a decomposition
$S_0(4)\cong\fsu(2)\oplus S^2V$, where $V=\CC^2=R(\pi_1)$ is
the natural module (the symmetric power is taken over $\CC$).
Moreover, $S^2R(\pi_1)\cong R(2\pi_1)$, and $\beta_{2\pi_1}=\RR$.
Thus $_\RR R(2\pi_1)=\fsu(2)$, and
$S_0(4)\cong\fsu(2)\oplus\fsu(2)\oplus\fsu(2)$.
Similarly as in Proposition \ref{ContainsSp(n)}, the
centralizer $\Sp(1)=\mathrm{Cen}_{\SO(4)}\SU(2)$ acts
transitively on the copies of $\fsu(2)$ in $S_0(4)$. Thus we have
only the cases 
\begin{align*}
\fh\cap S_0(4)\cong\fsu(2)&
\quad\text{ and }\quad
\fh\cong\fsl_2\CC \\
\fh\cap S_0(4)\cong\fsu(2)\oplus\fsu(2)&
\quad\text{ and }\quad
\fh\cong\fsp_4\RR \\
\fh\cap S_0(4)\cong S_0(4)&
\quad\text{ and }\quad
\fh=\fsl_4\RR.
\end{align*}
\qed
\end{Prop}
We summarize these results as follows.
\begin{Thm}
\label{TrsProj}
Let $G\leq\GL_m\RR$ be a closed subgroup, for $m\geq 3$.
Suppose that $G$ acts transitively on the projective space
$\RR\mathrm{P}^{m-1}$. Let $L=[G,G]$ denote the commutator
group of $G$. Then $G\leq\mathrm{Nor}_{\GL_m\RR}(L)$, and
there is a split short exact sequence
\[
\xymatrix{
1\ar[r] & L\ar[r] & \mathrm{Nor}_{\GL_m\RR}(L)\ar[r]^(.65){\leftarrow}
 & N\ar[r] & 1.}
\]
The following list gives all possibilities for
$L$ and the factor group $N=\mathrm{Nor}_{\GL_m\RR}(L)/L$.
We put $\RR_>=\{r\in\RR\mid r>0\}$.

\begin{center}
\begin{tabular}{l|ll}
$L$ & $m$ & $N$ \\ \hline
$\SO(2n)$        & $2n$      & $\RR_>$ \\
$\SO(2n+1)$      & $2n+1$    & $\RR^*$ \\
$\SU(n)$         & $2n$\quad($n\geq 3$) & $\CC^*$ \\
$\Sp(n)$         & $4n$      & $\SO(3)\cdot\RR_>$ \\
$\Sp(n)\cdot\Sp(1)$ & $4n$   & $\RR_>$ \\
$\mathrm{G}_2$   & $7$       & $\RR^*$ \\
$\Spin(7)$       & $8$       & $\RR_>$ \\
$\Spin(9)$       & $16$      & $\RR_>$ \\ \hline
$\SL_{2n}\RR$    & $2n$      & $\RR_>$ \\
$\SL_{2n+1}\RR$  & $2n+1$    & $\RR^*$ \\
$\Sp_{2n}\RR$    & $2n$      & $\RR_>$ \\
$\SL_n\CC$       & $2n$      & $\CC^*\rtimes\ZZ/2$ \\
$\Sp_{2n}\CC$    & $4n$      & $\CC^*\rtimes\ZZ/2$ \\
$\SL_n\HH$       & $4n$      & $\SO(3)\cdot\RR_>$ \\
$\SL_n\HH\cdot\Sp(1)$ & $4n$ & $\RR_>$ \\
$\Spin_{9,1}\RR$ & $16$      & $\RR_>$
\end{tabular}
\end{center}

\proof
The possibilities for the group $L$ were determined in Propositions
\ref{ContainsSO(n)}--\ref{ContainsSp(n)}. In each case, it is not
difficult to determine the normalizer and to construct a splitting of
the exact sequence.
\qed
\end{Thm}
The proof of the next lemma is straight-forward; it can be used to
derive a list of all connected transitive groups.
\begin{Lem}
\label{Spiral}
Let $S\leq\HH^*$ be a closed noncompact connected
1-dimensional subgroup. Up to conjugation, $S$ is of the form
\[
S=S_{a}=\{e^{t(1+\ti a)}\mid t\in\RR\}=
\{xe^{\ti a\mathrm{ln}(x)}\mid x\in\RR_{>}\}
\]
for some real number $a\in\RR$. The group is central in $\HH^*$
if and only if $a=0$, and $S_0=\RR_{>}$.

\proof
Let $\mathfrak{s}\leq\RR\oplus\fsu(2)$ denote the Lie algebra of
$S$. After conjugation with a suitable element $g\in\HH^*$, we may
assume that
$\mathfrak{s}\leq \RR\oplus\fso(2)\cong\RR\oplus\ti\RR$.
Let $(x,\ti y)\in\mathfrak{s}$
be a generator. Then $\exp(tx,\ti ty)=e^{tx}e^{\ti ty}$.
Since we assumed that $S$ is not compact, $x\neq 0$ and we may put
$(x,y)=(1,a)$.
\qed
\end{Lem}
In the 2-dimensional case, there are the following possibilities.
\begin{Lem}
\label{TrsR2}
Let $H\leq\SL_2\RR$ be a connected
group acting transitively on the nonzero vectors.
Then $H=\SL_2\RR$, or $H$ is conjugate to $\CC^*\leq\SL_2\RR$.

\proof
We have $3=\dim(\SL_2\RR)\geq\dim(H)\geq 2$, and $H$ is reductive. Thus
$H$ is abelian if $H\neq\SL_2\RR$. In the abelian case, $H$ acts
regularly and is thus homeomorphic to $\CC^*$; in particular, it contains
a torus $\SO(2)$. The connected centralizer of $\SO(2)$ is
$\CC^*=\mathrm{Cen}_{\SL_2\RR}\SO(2)$.
\qed
\end{Lem}
Combining the results of this section, we obtain the following
final result.
\begin{Thm}
\label{TrsLinGrp}
Let $H\leq\GL_m\RR$ be closed connected subgroup which acts
transitively on $\RR^m\setminus\{0\}$. Up to conjugation, $H$
is one of the groups listed in (a), (b), (c) below.

\smallskip
(a) $[H,H]$ is compact and $m\geq 3$.
\begin{center}
\begin{tabular}{l|lll}
$[H,H]$ & $m$ & $H$ \\ \hline
$\SO(n)$ & $n$ & $\SO(n)\cdot\RR_{>}$ \\
$\SU(n)$ & $2n$ & $\SU(n)\cdot S_a$,
\ $\SU(n)\cdot\CC^*$ & $(+)$ \\
$\Sp(n)$ & $4n$ & $\Sp(n)\cdot S_a$,
\ $\Sp(n)\cdot\CC^*$ & $(+)$ \\
$\Sp(n)\cdot\Sp(1)$ & $4n$ & $\Sp(n)\cdot\HH^*$ \\
$\mathrm{G}_2$ & $7$ & $\mathrm{G}_2\cdot\RR_{>}$ \\
$\Spin(7)$ & $8$ & $\Spin(7)\cdot\RR_{>}$ \\
$\Spin(9)$ & $16$ & $\Spin(9)\cdot\RR_{>}$
\end{tabular}
\end{center}
Here, $a$ can be any real number.

\smallskip
(b) $[H,H]$ is noncompact and $m\geq 3$.
\begin{center}
\begin{tabular}{l|lll}
$[H,H]$ & $m$ & $H$ \\ \hline
$\SL_n\RR$ & $n$ & $\SL_n\RR$,
\ $\SL_n\RR\cdot\RR_{>}$ \\
$\SL_n\CC$ & $2n$ & $\SL_n\CC$,
\ $\SL_n\CC\cdot\U(1)$,
\ $\SL_n\CC\cdot S_a$, 
\ $\GL_n\CC$ & $(+)$ \\
$\SL_n\HH$ & $4n$ & $\SL_n\HH$,
\ $\SL_n\HH\cdot\U(1)$,
\ $\SL_n\HH\cdot S_a$,
\ $\SL_n\HH\cdot\CC^*$ & $(+)$ \\
$\SL_n\HH\cdot\Sp(1)$ & $4n$ &
$\SL_n\HH\cdot\Sp(1)$,
\ $\SL_n\HH\cdot\HH^*$ \\
$\Sp_{2n}\RR$ & $2n$ & $\Sp_{2n}\RR$,
\ $\Sp_{2n}\RR\cdot\RR_{>}$ \\
$\Sp_{2n}\CC$ & $4n$ & $\Sp_{2n}\CC$,
\ $\Sp_{2n}\CC\cdot\U(1)$,
\ $\Sp_{2n}\CC\cdot S_a$,
\ $\Sp_{2n}\CC\cdot\CC^*$ & $(+)$ \\
$\Spin_{9,1}\RR$ & $16$ & $\Spin_{9,1}\RR$,
\ $\Spin_{9,1}\RR\cdot\RR_{>}$
\end{tabular}
\end{center}
Again, $a$ can be any real number.

\smallskip
(c) For $m=1,2$ there are only the following possibilites.
\begin{center}
\begin{tabular}{l|ll}
$m$ & $H$ \\ \hline
$1$ & $\RR^*$ \\
$2$ & $\CC^*$ & $(+)$ \\
& $\SL_2\RR$,\ $\SL_2\RR\cdot\RR_{>}$
\end{tabular}
\end{center}

\smallskip
(d) If $H$ preserves a complex structure on $\RR^m$, then
$H$ is one of the groups in (a), (b), (c) marked with $(+)$.

\proof
For $m\geq3$, the possibilities for the semisimple group
$L=[H,H]$ are determined in Theorem \ref{TrsProj}.
A direct inspection (combined with Lemma \ref{Spiral})
yields the list. For $m=2$ we use Lemma \ref{TrsR2}, and
the case $m=1$ is trivial. Finally, (d) follows by direct inspection
of the actions.
\qed
\end{Thm}
Using the results in this section, it is not difficult to obtain
a list of linear groups acting transitively on the point set
$\FF\mathrm{P}^k$ of a projective space over
$\FF=\CC,\HH$, see V\"olklein \cite{Volk} Satz 2.
Also, the possibilites for closed, but
not necessarily connected groups can be determined using
\ref{TrsProj}. We leave this to the reader.

\section{Locally compact Moufang sets}

A \emph{Moufang set} is a tripel $(G,U,X)$, where
$G$ is a 2-transitive permutation group acting on $X$,
and $U\unlhd G_x$ is
a normal subgroup of a stabilizer $G_x$ acting regularly on
$X\setminus\{x\}$. Moufang sets were introduced
by Tits in \cite{TitsDur}; they are also known as
\emph{split doubly transitive groups}. Note that the special case
$G_x=U$ is the same as a \emph{sharply 2-transitive group}.
See also Kramer \cite{CG} Section 1.8.

In this section, we determine all Moufang sets, where
$G$ is an effective 2-transitive Lie group and $U$ is closed in $G_x$.
We call such a Moufang set a \emph{locally compact and connected
Moufang set}.
According to the classification of 2-transitive Lie groups,
we distinguish three cases: the case where $G^\circ$ is simple
and of rank 1, the case where $G^\circ$ is simple and of higher
rank, and the affine case where $X=\RR^n$.

Recall from Section \ref{CompactCase}
the \emph{Iwasawa decomposition} of a simple Lie group
\[
H=KAU
\]
and the corresponding minimal parabolic
\[
B=K_0AU,
\]
where $K_0=\mathrm{Cen}_K(A)$ is the reductive anisotropic kernel.
\begin{Prop}
Let $G$ be a 2-transitive Lie group, with $G^\circ=H$ simple and
of rank $1$, as in Theorem \ref{CompactClass} (a).
As above, let $U$ denote the unipotent radical of a minimal
parabolic $B\subseteq H$. Then $(G,U,H/B)$ is a
locally compact connected Moufang set.

\proof
This is clear from the Iwasawa decomposition. Let $X=G/G_x=H/B$, and
let $y\in X\setminus\{x\}$. Then $H_{x,y}=B_y=K_0A$.
The group $U$ is normal in $B=H_x$, and intersects the two-point
stabilizer $K_0A$ trivially. Thus it acts regularly on $X\setminus\{x\}$.
\qed
\end{Prop}
If $G^\circ$ is a simple 2-transitive Lie
group of rank at least $2$, then there
is no way of making $X=G/G_x$ into a Moufang set. The following proof
was pointed out by Hendrik Van Maldeghem, replacing a more topological
(and more complicated) argument of mine.
Let $X$ be the point set of a desarguesian projective space of rank
at least $2$, or the point space of
a projective Moufang plane.
Let $G$ be a group
of automorphisms of the projective space, containing all elations.
Then $(G,X)$ cannot be made into a Moufang set.
To see this, let $p\in X$ and assume that $U\unlhd G_p$ is
a normal subgroup acting regularly on $X\setminus\{p\}$.
Let $u\in U$, and let $\tau$ be an elation with center $p$. Then
$u\tau u^{-1}$ is also an elation with center $p$, and so is
the commutator $[u,\tau]$.
If we choose $u$, $\tau$ in such a way that $u$ does not fix
the axis of $\tau$ (which is possible, since the rank of the projective
space is at least $2$), then
$[u,\tau]\in U$ is a nontrivial elation with center $p$.
Since $U$ is normal, $U$ contains all elations with center $p$.
These elations form an abelian
normal subgroup of $G_p$ which is, however, not regular on $X$.
\begin{Prop}
None of the groups in Theorem \ref{CompactClass} (b)
can be made into a Moufang set.
\qed
\end{Prop}
Finally, we consider the question of uniqueness in the case where
$H=G^\circ$ is simple and of rank $1$. The question is thus
if $G_x$ admits a regular normal subgroup $V$ different from $U$.
Since $X\setminus\{x\}$ is connected, we have $V\unlhd B=H_x=K_0AU$.
So $V\cap K_0A=1$ (by regularity of $V$) and $K_0AV=B$ (by
transitivity). Let $\mathfrak u$ denote the Lie algebra of $U$,
and $\mathfrak v$ the Lie algebra of $V$. We have to prove that
$\mathfrak v=\mathfrak u$.
We decompose the Lie algebra $\mathfrak b$ of $B$ into irreducible
$K_0A$-modules. As a $K_0A$-module, $\mathfrak k_0\oplus\mathfrak a$
decomposes as
$[\mathfrak k_0,\mathfrak k_0]\oplus\mathrm{Cen}(\mathfrak k_0)
\oplus\mathfrak a$.
Direct inspection shows that no $K_0A$-submodule of $\mathfrak u$
is isomorphic to a $K_0A$-submodule of $\mathfrak k_0\oplus\mathfrak a$.
This proves uniqueness of $\mathfrak u$.
\begin{Lem}
Each of the groups in Theorem \ref{CompactClass} (a)
is in a unique way a
locally compact and connected Moufang set, i.e. the
data $(G,X)$ determine the closed regular normal subgroup
$U\unlhd G_x$ uniquely.
\qed
\end{Lem}
It remains to consider the case where $X$ is noncompact.
Then $G$ is a semidirect product $G=G_x\ltimes\RR^m$, and
$U\unlhd G_x$ is a normal subgroup acting regularly on the
nonzero vectors in $\RR^m$.
It follows that the chosen maximal compact subgroup $K\subseteq U$
acts regularly on the sphere $\SS^{m-1}$. Direct inspection of the
list in Theorem \ref{TrsSphere}
shows that this happens only for $m=1,2,4$, and
$[U,U]$ is one of the groups $1,\SO(2),\Sp(3)$.
We obtain the following result, which is originally due to
Kalscheuer \cite{Kal} and was re-proved by Tits \cite{TitsKal}
and Grundh\"ofer \cite{GrundKal}.
\begin{Thm}
\label{Sh2Trs}
Let $U\leq\GL_m\RR$ be a closed subgroup which acts regularly 
(i.e. sharply transitively) on
$\RR^m\setminus\{0\}$. Then $m=1,2,4$, and $U$ is one of the
groups $\RR^*$, $\CC^*$, $\HH^*=\Sp(1)\cdot S_0$, or
$\Sp(1)\cdot S_a$, for $a\neq 0$.

\proof
The result follows by direct inspection of the tables in 
Theorem \ref{TrsLinGrp}.
\qed
\end{Thm}
Since $U\unlhd H$ is a normal subgroup, it remains to determine
the normalizer $\mathrm{Nor}_{\GL_m\RR}(U)$. The following table
shows the corresponding quotients $\mathrm{Nor}_{\GL_m\RR}(U)/U$.
\begin{center}
\begin{tabular}{l|ll}
$U$ & $m$ & $\mathrm{Nor}_{\GL_m\RR}(U)/U$ \\ \hline
$\RR^*$        & $1$      & $1$ \\
$\CC^*$        & $2$      & $\ZZ/2$ \\
$\HH^*$        & $4$      & $\SO(3)\cdot\RR_>$ \\
$\Sp(1)\cdot S_a$ & $4$   & $\RR_>$ \quad($a\neq 0$)
\end{tabular}
\end{center}
Using this, it is not difficult to determine the locally compact
and connected Moufang sets for noncompact $X$.
Also, it is easy to see that the pair $(G,X)$ determines $U$ for
$m=1,2$. This is \emph{not} true for $m=4$: if $H=\Sp(1)\cdot \CC^*$
is the group consisting of all maps $x\mapstoo hxc$, for
$h\in\Sp(1)$ and $c\in\CC^*$, then $\Sp(1)\cdot S_a=U\unlhd H$ is
a regular normal subgroup for any choice of $a\in\RR$. 

\medskip
\textbf{Acknowledgements.}
I would like to thank Theo Grundh\"ofer, Rainer L\"owen, Bernhard
M\"uhl\-herr,
Katrin Tent, Hendrik Van Maldeghem and the referee for helpful remarks, and
Jacques Tits for some comments on the history of the whole subject.

\small
\makeatletter
\parbox[t]{10cm}{\begin{raggedright}
Linus Kramer \\
Mathematisches Institut \\
Universit\"at W\"urzburg \\
Am Hubland \\
D--97074 W\"urzburg \\
Germany \\
\footnotesize{\tt kramer@mathematik.uni-wuerzburg.de}\
\end{raggedright}}

\end{document}